April 15, 1998

# Ergodicity of Hard Spheres in a Box

Nándor Simányi[1]
Bolyai Institute of Mathematics
The University of Szeged
Szeged, Aradi Vértanuk tere 1. H-6720 Hungary
E-mail: simanyi@math.u-szeged.hu

Abstract. We prove that the system of two hard balls in a $\nu$-dimensional ($\nu \geq 2$) rectangular box is ergodic and, therefore, actually it is a Bernoulli flow.

## 1. Introduction

### The Model and the Theorem

Let us consider the billiard system of two hard balls with unit mass and radius $r$ ($0 < r < 1/4$) moving uniformly in the $\nu$-dimensional ($\nu \geq 2$) Euclidean container

$$\mathbf{C} = [-r,\, 1+r]^k \times \mathbb{T}^{\nu-k} = [-r,\, 1+r]^k \times \left(\mathbb{R}^{\nu-k}/\mathbb{Z}^{\nu-k}\right)$$

($0 \leq k \leq \nu$) and bouncing back elastically at each other and at the boundary $\partial \mathbf{C}$ of $\mathbf{C}$. Denote the center of the $i$-th ball ($i = 1, 2$) by $q_i$, and its time derivative by $v_i = \dot{q}_i$. Also denote by $\mathcal{A} = \{1, 2, \ldots, k\}$ the set of the first $k$ (i. e. the non-periodic) coordinate axes of $\mathbf{C}$, and by $\pi_2(.)$ the projection (of a position or a velocity) into the second, periodic factor $\mathbb{T}^{\nu-k}$ of the container $\mathbf{C}$. We use the usual reductions $2E = ||v_1||^2 + ||v_2||^2 = 1$, $\pi_2(q_1 + q_2) = \pi_2(v_1 + v_2) = 0$. Plainly, the configuration space is the set

$$\mathbf{Q} = \big\{(q_1,\, q_2) \in \left([0,\, 1]^k \times \mathbb{T}^{\nu-k}\right) \times \left([0,\, 1]^k \times \mathbb{T}^{\nu-k}\right) \,\big|$$
$$\text{dist}(q_1,\, q_2) \geq 2r \text{ and } \pi_2(q_1 + q_2) = 0\big\}$$

with a connected interior. The phase space $\mathbf{M}$ of the arising semi-dispersive billiard flow is almost the unit tangent bundle of $\mathbf{Q}$. The only modification to the unit tangent bundle is that the incoming and outgoing velocities at $\partial \mathbf{Q}$ are glued together, just as it is prescribed by the law of elastic collisions: the post-collision (outgoing) velocity is the reflected image of the pre-collision (incoming) velocity across the tangent hyperplane of $\partial \mathbf{Q}$ at the point of collision. We call the arising semi-dispersive billiard flow *the standard $(\nu, k, r)$ model*, or *the standard $(\nu, k, r)$ flow*.

---

[1]Research supported by the Hungarian National Foundation for Scientific Research, grants

**Theorem.** *For every triple $(\nu, k, r)$ ($\nu \geq 2$, $0 \leq k \leq \nu$, $0 < r < 1/4$) the standard $(\nu, k, r)$ model is ergodic, hence it is actually a Bernoulli flow, see [3], [5], or [15].*

**Remark.** It becomes clear from the upcoming proof of this theorem that the assumptions

(1) on the equality of the side lengths of the container **C**,
(2) on the equality of the masses of balls, and
(3) on the equality of the radii of balls

are not essential, but merely notational simplifications. The proof easily carries over to the general case when these equalities do not hold.

## Some Historical Notes

Hard ball systems or, a bit more generally, mathematical billiards constitute an important and quite interesting family of dynamical systems being intensively studied by dynamicists and researchers of mathematical physics, as well. These dynamical systems pose many challenging mathematical questions, most of them concerning the ergodic (mixing) properties of such systems. The introduction of hard ball systems and the first major steps in their investigations date back to the 40's and 60's, see Krylov's paper [13] and Sinai's groundbreaking works [20], [21]. In the articles [21] and [2] Bunimovich and Sinai prove the ergodicity of two hard disks in the two-dimensional unit torus $\mathbb{T}^2$. The generalization of this result to higher dimensions $\nu > 2$ took fourteen years, and was done by Chernov and Sinai in [22]. Although the model of two hard balls in $\mathbb{T}^\nu$ is already rather involved technically, it is still a so called strictly dispersive billiard system, i. e. such that the smooth components of the boundary $\partial \mathbf{Q}$ of the configuration space are strictly convex from inside **Q**. The billiard systems of more than two hard spheres in $\mathbb{T}^\nu$ are no longer strictly dispersive, but just dispersive (strict convexity of the smooth components of $\partial \mathbf{Q}$ is lost, merely convexity persists!), and this circumstance causes a lot of additional technical troubles in their study. In the series of my joint papers with A. Krámli and D. Szász [9–12] we developed several new methods, and proved the ergodicity of more and more complicated semi-dispersive billiards culminating in the proof of the ergodicity of four billiard balls in the torus $\mathbb{T}^\nu$ ($\nu \geq 3$), [12]. Then, in 1992, Bunimovich, Liverani, Pellegrinotti and Sukhov [1] were able to prove the ergodicity for some systems with an arbirarily large number of hard balls. The shortcoming of their model, however, is that, on one hand, they restrict the types of all feasible ball-to-ball collisions, on the other hand, they introduce some extra scattering effect with the collisions at the strictly convex wall of the container. The only result with an arbitrarily large number of spheres in a flat unit torus $\mathbb{T}^\nu$ was achieved in [16–17], where I managed to prove the ergodicity (actually, the K-mixing property) of $N$ hard balls in $\mathbb{T}^\nu$, provided that $N \leq \nu$. The annoying shortcoming of that result is that the larger the number of balls $N$ is, larger and larger dimension $\nu$ of the ambient container is required by the method of the proof.

On the other hand, if someone considers a hard sphere system in an elongated torus which is long in one direction but narrow in the others, so that the spheres must keep their cyclic order in the "long direction" (Sinai's "pencase" model), then spheres are now restricted to neighbouring pairs. The hyperbolicity of such models in three dimensions and the ergodicity in dimension four have been proved in [18].

The positivity of the metric entropy for several systems of hard spheres can be proven easily, as was shown in the paper [26]. The articles [14] and [27] are nice surveys describing a general setup leading to the technical problems treated in the present paper. For a comprehensive survey of the results and open problems in this field, see [24].

Finally, in our latest joint venture with D. Szász [19] we prevailed over the difficulty caused by the low value of the dimension $\nu$ by developing a brand new algebraic approach for the study of hard ball systems. That result, however, only establishes hyperbolicity (nonzero Lyapunov exponents almost everywhere) for $N$ balls in $\mathbb{T}^\nu$. The ergodicity is a bit longer shot.

None of the above results took up the problem of handling hard balls in physically more realistic containers, e. g. rectangular boxes. The extra technical hardship in their investigation is caused by the loss of the total momentum and center of mass. This amounts to the increase in the dimension of the configuration (phase) space without any additional scattering effect as a compensation. The problem of proving ergodicity for $N$ hard spheres in a $\nu$-dimensional rectangular box is so difficult that we were only able to achieve that goal in the case $N = 2$.

## The Strategy of the Proof

After reviewing the necessary technical skills in Section 2, in the subsequent section we introduce the concept of combinatorial richness of the symbolic collision structure of a trajectory (segment), and then we prove that — apart from some codimension-two, smooth submanifolds of the phase space — such a combinatorial richness implies the sufficiency (hyperbolicity) of the trajectory.

The proof of the theorem goes on by an induction on the number $k = 0, 1, \ldots, \nu$. Section 4 uses the inductive hypothesis (the statement of the theorem for $k-1$) and proves that the set of combinatorially non-rich phase points is *slim*, i. e. it can be covered by a countable family of closed, zero measure sets with codimension at least two.

Section 5 contains the inductive proof of the theorem based upon the preceding two sections and the celebrated Theorem on Local Ergodicity for semi-dispersive billiards by Chernov and Sinai, [22].

The closing Appendix is a brief overview of a special orthogonal cylindric billiard which emerged in Section 4.

## 2. Prerequisites

### Semi-dispersive Billiards

A billiard is a dynamical system describing the motion of a point particle in a connected, compact domain $\mathbf{Q} \subset \mathbb{R}^d$ or $\mathbf{Q} \subset \mathbb{T}^d = \mathbb{R}^d/\mathbb{Z}^d$, $d \geq 2$, with a piecewise $C^2$-smooth boundary. Inside $\mathbf{Q}$ the motion is uniform, whereas the reflection at the boundary $\partial \mathbf{Q}$ is elastic (the angle of reflection equals the angle of incidence). Since of our system can be identified with the unit tangent bundle over $\mathbf{Q}$. Namely, the configuration space is $\mathbf{Q}$, while the phase space is $\mathbf{M} = \mathbf{Q} \times \mathbb{S}^{d-1}$, where $\mathbb{S}^{d-1}$ is the unit $d-1$-sphere. In other words, every phase point $x$ is of the form $(q, v)$ where $q \in \mathbf{Q}$ and $v \in \mathbb{S}^{d-1}$ is a tangent vector at the footpoint $q$. The natural projections $\pi : \mathbf{M} \to \mathbf{Q}$ and $p : \mathbf{M} \to \mathbb{S}^{d-1}$ are defined by $\pi(q, v) = q$ and $p(q, v) = v$, respectively.

Suppose that $\partial \mathbf{Q} = \cup_1^k \partial \mathbf{Q}_i$, where $\partial \mathbf{Q}_i$ are the smooth components of the boundary. Denote $\partial \mathbf{M} = \partial \mathbf{Q} \times \mathbb{S}^{d-1}$, and let $n(q)$ be the unit normal vector of the boundary component $\partial \mathbf{Q}_i$ at $q \in \partial \mathbf{Q}_i$ directed inwards $\mathbf{Q}$.

The flow $\{S^t : t \in \mathbb{R}\}$ is determined for the subset $\mathbf{M}' \subset \mathbf{M}$ of phase points whose trajectories never cross the intersections of the smooth pieces of $\partial \mathbf{Q}$ and do not contain an infinite number of reflections in a finite time interval. If $\mu$ denotes the (normalized) Liouville measure on $\mathbf{M}$, i.e. $d\mu(q, v) = \text{const} \cdot dq \cdot dv$, where $dq$ and $dv$ are the differentials of the Lebesgue measures on $\mathbf{Q}$ and on $\mathbb{S}^{d-1}$, respectively, then under certain conditions $\mu(\mathbf{M}') = 1$ and $\mu$ is invariant [8]. The interior of the phase space $\mathbf{M}$ can be endowed with the natural Riemannian metric.

The dynamical system $(\mathbf{M}, \{S^t\}, \mu)$ is called a *billiard*. Notice, that $(\mathbf{M}, \{S^t\}, \mu)$ is neither everywhere defined nor smooth.

The main object of the present paper is a particularly interesting class of billiards: the *semi-dispersive billiards* where, for every $q \in \partial \mathbf{Q}$ the second fundamental form $K(q)$ of the boundary is positive semi-definite. If, moreover, for every $q \in \partial \mathbf{Q}$ the second fundamental form $K(q)$ is positive definite, then the billiard is called a *dispersive billiard*.

As it is pointed out in previous works on billiards, the dynamics can only be defined for trajectories where the moments of collisions do not accumulate in any finite time interval (cf. Condition 2.1 of [10]). An important consequence of Theorem 5.3 of [25] is that – for semi-dispersing billiards – there are *no trajectories at all with a finite accumulation point of collision moments*, see also [6].

### Convex Orthogonal Manifolds

In the construction of invariant manifolds a crucial role is played by the time evolution equation for the second fundamental form of codimension-one submanifolds in $\mathbf{Q}$ orthogonal to the velocity component $p(x)$ of a phase point $x$. Let $x = (q, v) \in \mathbf{M} \setminus \partial \mathbf{M}$, and consider a $C^2$-smooth, codimension-one submanifold $\tilde{\mathcal{O}} \subset \mathbf{Q} \setminus \partial \mathbf{Q}$ such that $q \in \tilde{\mathcal{O}}$ and $v = p(x)$ is a unit normal vector to $\tilde{\mathcal{O}}$ at $q$. Denote by $\mathcal{O}$ the normal section of the unit tangent bundle of $\mathbf{Q}$ restricted to $\tilde{\mathcal{O}}$. (The manifold $\mathcal{O}$ is uniqely defined by the orientation $(q, v) \in \mathcal{O}$.) We call $\mathcal{O}$ a *local orthogonal manifold* with support $\tilde{\mathcal{O}}$.

Recall that *the second fundamental form* $B_{\mathcal{O}}(x)$ of $\mathcal{O}$ (or $\tilde{\mathcal{O}}$) at $x$ is defined through
$$v(q + \delta q) - v(q) = B_{\mathcal{O}}(x) \cdot \delta q + o\left(\|\delta q\|\right),$$
and it is a self-adjoint operator acting in the $(d-1)$-dimensional tangent hyperplane of $\tilde{\mathcal{O}}$ at $q$.

A local orthogonal manifold $\mathcal{O}$ is called *convex* if $B_{\mathcal{O}}(y) \geq 0$ for every $y \in \mathcal{O}$.

Recall that the common tangent space $\mathcal{T}_q \mathbf{Q}$ of the parallelizable configuration

## Neutral Vectors, Advance and Sufficiency

Consider a *non-singular* trajectory segment $S^{[a,b]}x$ in a semi-dispersive billiard. Suppose that $a$ and $b$ are *not moments of collision*.

**Definition 2.1.** *The neutral space $\mathcal{N}_0(S^{[a,b]}x)$ of the trajectory segment $S^{[a,b]}x$ at time zero ($a \leq 0 \leq b$) is defined by the following formula:*

$$\mathcal{N}_0(S^{[a,b]}x) = \left\{ w = \left(\delta q_1^0, \delta q_2^0; \delta v_1^0, \delta v_2^0\right) \in \mathcal{T}_x \mathbf{M} \middle| \, \delta v_i^t = 0 \text{ for } i = 1, 2; \, t \in [a, b] \right\},$$

*where $w^t = (\delta q_1^t, \delta q_2^t; \delta v_1^t, \delta v_2^t) = (DS^t)(w) \in \mathcal{T}_{x^t} \mathbf{M}$ denotes the image of the tangent vector $w = w^0 \in \mathcal{T}_x \mathbf{M}$ under the linearization $DS^t$ of the map $S^t$, $t \in [a, b]$.*

**Remark–definition 2.2.** It follows easily from the above definition that for every tangent vector $w \in \mathcal{N}_0(S^{[a,b]}x)$ and for every moment $t \in [a, b]$ of a ball-to-ball collision we have

$$\delta q_1^{t-0} - \delta q_2^{t-0} = \alpha(w) \cdot \left[v_1^{t-0} - v_2^{t-0}\right] \text{ and}$$
$$\delta q_1^{t+0} - \delta q_2^{t+0} = \alpha(w) \cdot \left[v_1^{t+0} - v_2^{t+0}\right],$$

where $\alpha: \mathcal{N}_0(S^{[a,b]}x) \to \mathbb{R}$ is a linear functional, see also Section 2 of [10]. The functional $\alpha$ is called the *advance* of the collision $\sigma$ taking place at time $t$.

It is now time to bring up the basic notion of *sufficiency* of a trajectory (segment). This is the utmost important necessary condition for the proof of the Theorem on Local Ergodicity for semi-dispersive billiards, see Condition (ii) of Theorem 3.6 and Definition 2.12 in [10].

**Definition 2.3.**
(1) *The non-singular trajectory segment $S^{[a,b]}x$ ($a$ and $b$ are supposed not to be moments of collision) is said to be sufficient if and only if the dimension of $\mathcal{N}_t(S^{[a,b]}x)$ ($t \in [a, b]$) is minimal, i.e. $\dim\left(\mathcal{N}_t(S^{[a,b]}x)\right) = 1$.*
(2) *The trajectory segment $S^{[a,b]}x$ containing exactly one singularity is said to be sufficient if and only if both branches of this trajectory segment are sufficient.*

The above definition uses the notion of the

## Trajectory Branches

We are going to briefly describe the discontinuity of the semi-dispersive billiard flow $\{S^t\}$ caused by a collision with several intersecting, smooth components of the boundary $\partial \mathbf{Q}$ at time $t_0$. Assume first that the pre-collision velocity $v$ is given. What can we say about the possible post-collision velocity? Let us perturb the pre-collision phase point (at time $t_0 - 0$) infinitesimally, so that the collisions at $\sim t_0$ (with the several smooth components of the boundary $\partial \mathbf{Q}$) occur at infinitesimally different moments. By applying the collision laws to the arising finite sequence of collisions, we see that the post-collision velocity is fully determined by the time-ordering of the considered collisions. Therefore, the collection of all possible time-orderings of these collisions gives rise to a finite family of continuations clear that similar statements can be said regarding the evolution of a trajectory through a multiple collision *in reverse time*. Furthermore, it is also obvious that for any given phase point $x_0 \in \mathbf{M}$ there are two, $\omega$-high trees $\mathcal{T}_+$ and $\mathcal{T}_-$ such that $\mathcal{T}_+$ ($\mathcal{T}_-$) describes all possible continuations of the positive (negative) trajectory $S^{[0,\infty)}x_0$ ($S^{(-\infty,0]}x_0$). (For the definition of trees and for some of their applications to billiards cf. the beginning of Section 5 in [12].) It is also clear that all the possible continuations (branches) of the whole trajectory $S^{(-\infty,\infty)}x_0$ can be uniquely described by all possible pairs $(B_-, B_+)$ of $\omega$-high branches of the trees $\mathcal{T}_-$ and $\mathcal{T}_+$ ($B_- \subset \mathcal{T}_-, B_+ \subset \mathcal{T}_+$).

### Slim (negligible) Sets

We are going to summarize the basic properties of codimension-two subsets $A$ of a smooth manifold $\mathbf{M}$. Since these subsets $A$ are just those absolutely negligible in our dynamical discussions, we shall call them *slim*. As to a broader exposition of the issues, see [4] or Section 2 of [11].

Note that the dimension $\dim A$ of a separable metric space $A$ is one of the three classical notions of dimension: the covering, the small inductive, or the large inductive dimension. As it is known from general topology, all of them are the same for separable metric spaces.

**Definition 2.4.** *A subset $A$ of $\mathbf{M}$ is called slim if $A$ can be covered by a countable family of codimension-two (i. e. at least two), closed sets of $\mu$-measure zero, where $\mu$ is some smooth measure on $\mathbf{M}$. (See also Definition 2.12 of [11].)*

**Property 2.5.** *The collection of all slim subsets of $\mathbf{M}$ is a $\sigma$-ideal, that is, countable unions of slim sets and arbitrary subsets of slim sets are also slim.*

**Lemma 2.6.** *A subset $A \subset \mathbf{M}$ is slim if and only if for every $x \in A$ there exists an open neighborhood $U$ of $x$ in $\mathbf{M}$ such that $U \cap A$ is slim. (Locality, cf. Lemma 2.14 of [11].)*

**Property 2.7.** *A closed subset $A \subset \mathbf{M}$ is slim if and only if $\mu(A) = 0$ and $\dim A \leq \dim \mathbf{M} - 2$.*

**Lemma 2.8.** *If $A \subset M_1 \times M_2$ is a closed subset of the product of two manifolds, and for every $x \in M_1$ the set*

$$A_x = \{y \in M_2 : (x,y) \in A\}$$

*is slim in $M_2$, than $A$ is slim in $M_1 \times M_2$.*

The following lemmas characterize codimension-one and codimension-two sets.

**Lemma 2.9.** *For any closed subset $S \subset \mathbf{M}$ the following three conditions are equivalent:*
  (i) $\dim S \leq \dim \mathbf{M} - 2$;
  (ii) $\text{int} S = \emptyset$ *and for every open and connected set $G \subset \mathbf{M}$ the difference set $G \setminus S$ is also connected;*
  (iii) $\text{int} S = \emptyset$ *and for every point $x \in \mathbf{M}$ and for any neighborhood $V$ of $x$ in $\mathbf{M}$ there exists a smaller neighborhood $W \subset V$ of the point $x$ such that, for every pair of points $y, z \in W \setminus S$, there is a continuous curve $\gamma$ in the set*

(See Theorem 1.8.13 and Problem 1.8.E of [4].)

**Lemma 2.10.** *For any subset $S \subset \mathbf{M}$ the condition $\dim S \leq \dim \mathbf{M} - 1$ is equivalent to $\text{int} S = \emptyset$. (See Theorem 1.8.10 of [4].)*

We recall an elementary, but important result, Lemma 4.15 of [11]. Let $R_2$ be the set of phase points $x \in \mathbf{M} \setminus \partial \mathbf{M}$ for which the trajectory $S^{(-\infty,\infty)}x$ has at least two singularities.

**Lemma 2.11.** *The set $R_2$ is a countable union of codimension-two (i. e. at least two), smooth submanifolds of $\mathbf{M}$.*

The last lemma in this section establishes the most important property of slim sets, which gives us the fundamental geometric tool to connect the open ergodic components of billiard flows.

**Lemma 2.12.** *If $\mathbf{M}$ is connected, then the complement $\mathbf{M} \setminus A$ of a slim set $A \subset \mathbf{M}$ necessarily contains an arcwise connected, $G_\delta$ set of full measure. (See Property 3 of Section 4.1 of [9]. The $G_\delta$ sets are, by definition, the countable intersections of open sets.)*

## 3. Richness Implies Sufficiency

### (The Case $k > 0$)

In this section we will be studying the neutral spaces and sufficiency of *non-singular* trajectory segments $\omega = S^{[a,b]}x_0$, where $a$ and $b$ are not moments of collision. We are interested in answering the following fundamental question: How does the neutral space $\mathcal{N}(\omega)$ depend on the symbolic collision sequence $\Sigma(\omega) = \Sigma$ of the given orbit segment $\omega = S^{[a,b]}x_0$, and what kind of combinatorial *richness* of $\Sigma$ guarantees sufficiency, i. e. that $\dim \mathcal{N}(\omega) = 1$?

### The Symbolic Sequence $\Sigma(\omega)$

Let us symbolically denote the really dispersive (i. e. ball-to-ball) collisions $\{1,2\}$ of the orbit segment $\omega$ by $\sigma_0, \sigma_1, \ldots, \sigma_n$ just as they follow each other in time, and denote by $a < t_0 = t(\sigma_0) < \cdots < t_n = t(\sigma_n) < b$ the moments when these collisions take place. For an arbitrary triple $(\beta, i, j)$ ($\beta = 1, 2$; $i = 1, \ldots, n$; $j = 1, \ldots, k$) we introduce the nonnegative integer $r(\beta, i, j)$ as the total number of collisions of the $\beta$-th ball with the flat boundary components $(q)_j = 0$ and $(q)_j = 1$ during the time interval $(t_{i-1}, t_i)$, i. e.

$$r(\beta, i, j) = \# \left\{ t \mid t_{i-1} < t < t_i \text{ and } (q_\beta^t)_j \cdot [1 - (q_\beta^t)_j] = 0 \right\}, \tag{3.1}$$

where the subscript $(.)_j$ denotes the $j$-th component of a vector. Then we define the following three sets of axes $Z_i(1), Z_i(2), Z_i \subset \mathcal{A}$ ($i = 1, \ldots, n$):

where $\triangle$ denotes the symmetric difference of sets. Recall from the introduction that the set of axes $\mathcal{A}$ is just the index set $\{1, 2, \ldots, k\}$, $(1 \leq k \leq \nu)$. For a set $Z \subset \{1, \ldots, \nu\}$ we denote by $R_Z : \mathbb{R}^\nu \to \mathbb{R}^\nu$ the orthogonal reflection across the subspace spanned by the coordinate axes $j \in \{1, \ldots, \nu\} \setminus Z$. It is then clear that

$$v_\beta^{t_i-0} = R_{Z_i(\beta)} v_\beta^{t_{i-1}+0} \tag{3.3}$$

for $i = 1, \ldots, n$, $\beta = 1, 2$. We call the sequence

$$\Sigma = \Sigma(\omega) = (\sigma_0, Z_1, \sigma_1, Z_2, \sigma_2, \ldots, Z_n, \sigma_n)$$

the symbolic collision sequence.

**Remark.** We note here that in the notation of $\Sigma$ the presence of $\sigma_i$'s is obviously redundant. The reason for still including them was simply to indicate the existence of the dispersive collisions between the non-dispersive ball-to-wall collisions. Secondly, the above notion of symbolic collision sequence differs from the usual concept where one encodes the itinerary of the orbit segment by showing the sequence of different smooth components of the boundary of the configuration space at which the collisions took place. This new concept of the symbolic collision sequence is a reduced version of the traditional one, containing less information on the orbit segment.

In the actual computations of the neutral space $\mathcal{N} = \mathcal{N}(\omega)$ not all sets $Z_i(\beta)$, but only their symmetric differences

$$Z_i = Z_i(1) \triangle Z_i(2)$$

will play a role, $i = 1, \ldots, n$. We denote by $\alpha_i = \alpha(\sigma_i) : \mathcal{N}(\omega) \to \mathbb{R}$ $(i = 0, \ldots, n)$ the advance functional corresponding to the collision $\sigma_i$, see also Section 2.

**Definition 3.4.** We say that the symbolic collision sequence

$$\Sigma = (\sigma_0, Z_1, \sigma_1, \ldots, Z_n, \sigma_n)$$

is *rich* (or, *combinatorially rich*) if
  (i) $\bigcup_{i=1}^n Z_i = \mathcal{A}$ and
  (ii) there exists an index $i$ for which $0 < |Z_i| < \nu$.
(We note that, since always $|Z_i| \leq k$, the second property actually follows from the first one if $k < \nu$.)

**Key Lemma 3.5.** *Let $\Sigma$ be a combinatorially rich symbolic sequence. Then there exists a countable union $E = E(\Sigma) \subset \mathbf{M}$ of proper, smooth submanifolds of the phase space $\mathbf{M}$ with the following property:*
   *If $\Sigma(\omega) = \Sigma$ for a non-singular orbit segment $\omega = S^{[a,b]} x_0$ and $x_0 \notin E$, then the orbit segment $\omega$ is sufficient, i. e. $\dim \mathcal{N}(\omega) = 1$.*

**Remark.** It will be clear from the proof of the key lemma that, whenever the exceptional phase points form a codimension-one submanifold $J \subset \mathbf{M}$, then such a manifold $J$ must be defineable by an equation $(v_1^t)_j - (v_2^t)_j = 0$ with some given $j \in \{1, \ldots, \nu\}$ and $t \in (a, b)$.

Most of the rest of this section will be devoted to the proof of the key lemma. The argument will be split into several lemmas.

**Lemma 3.6.** *Assume that $n = 1$ and $|Z_1| < \nu$. Consider an arbitrary neutral vector $w = \left(\delta q_1^{t_0+0}, \delta q_2^{t_0+0}\right) \in \mathcal{N}(\omega)$ with advances $\alpha_0(w)$, $\alpha_1(w)$. We claim that $\alpha_0(w) = \alpha_1(w)$, and*

$$P_{Z_1}\left(\delta q_\beta^{t_0+0}\right) = \alpha_0(w) \cdot P_{Z_1}\left(v_\beta^{t_0+0}\right), \quad \beta = 1, 2,$$

*unless $P_{\overline{Z}_1}\left(v_1^{t_0+0} - v_2^{t_0+0}\right) = 0$. Here $P_{Z_1} = \dfrac{1}{2}(I - R_{Z_1})$ is the orthogonal projection onto the subspace spanned by the axes in $Z_1$ and $\overline{Z}_1 = \{1, 2, \ldots, \nu\} \setminus Z_1$ is the complement of $Z_1$.*

**Proof.** The neutrality of $w$ means that

$$\delta q_1^{t_0+0} - \delta q_2^{t_0+0} = \alpha_0(w) \cdot \left[v_1^{t_0+0} - v_2^{t_0+0}\right] \quad \text{and}$$
$$\delta q_1^{t_1-0} - \delta q_2^{t_1-0} = \alpha_1(w) \cdot \left[v_1^{t_1-0} - v_2^{t_1-0}\right],$$

see 2.2. By using these equations and the transformation equation (3.3) (for the velocities and the position variations $\delta q_\beta$, as well) we get

$$\alpha_1(w) \cdot P_{\overline{Z}_1}\left(v_1^{t_1-0} - v_2^{t_1-0}\right) = P_{\overline{Z}_1}\left(\delta q_1^{t_1-0} - \delta q_2^{t_1-0}\right) = R_{Z_1(1)} P_{\overline{Z}_1}\left(\delta q_1^{t_0+0} - \delta q_2^{t_0+0}\right)$$
$$= \alpha_0(w) \cdot R_{Z_1(1)} P_{\overline{Z}_1}\left(v_1^{t_0+0} - v_2^{t_0+0}\right) = \alpha_0(w) \cdot P_{\overline{Z}_1}\left(v_1^{t_1-0} - v_2^{t_1-0}\right).$$

It follows immediately that $\alpha_0(w) = \alpha_1(w)$, unless the vector

$$P_{\overline{Z}_1}\left(v_1^{t_1-0} - v_2^{t_1-0}\right) = R_{Z_1(1)} P_{\overline{Z}_1}\left(v_1^{t_0+0} - v_2^{t_0+0}\right)$$

is zero, i. e. $P_{\overline{Z}_1}\left(v_1^{t_0+0} - v_2^{t_0+0}\right) = 0$.

Assume now that $\alpha_0(w) = \alpha_1(w)$. (Which is the typical situation, as we have seen.) Again, by using the neutrality equations (the first two equations in this proof) and (3.3), we obtain

$$P_{Z_1}\left(\delta q_1^{t_1-0} + \delta q_2^{t_1-0}\right) = R_{Z_1(1)} P_{Z_1} \delta q_1^{t_0+0} + R_{Z_1(2)} P_{Z_1} \delta q_2^{t_0+0} =$$
$$= R_{Z_1(1)} P_{Z_1} \delta q_1^{t_0+0} - R_{Z_1(1)} P_{Z_1} \delta q_2^{t_0+0} = \alpha_0(w) \cdot R_{Z_1(1)} P_{Z_1}\left(v_1^{t_0+0} - v_2^{t_0+0}\right) =$$
$$= \alpha_0(w) \cdot \left(P_{Z_1} v_1^{t_1-0} + P_{Z_1} v_2^{t_1-0}\right) = P_{Z_1}\left[\alpha_0(w) \cdot \left(v_1^{t_1-0} + v_2^{t_1-0}\right)\right].$$

On the other hand,

$$P_{Z_1}\left(\delta q_1^{t_1-0} - \delta q_2^{t_1-0}\right) = P_{Z_1}\left[\alpha_1(w) \cdot \left(v_1^{t_1-0} - v_2^{t_1-0}\right)\right].$$

The last two equations and $\alpha_0(w) = \alpha_1(w)$ together yield

$$P_{Z_1} \delta q_\beta^{t_1-0} = \alpha_0(w) \cdot P_{Z_1} v_\beta^{t_1-0}, \quad \beta = 1, 2.$$

Now applying the reflection $R_{Z_1(\beta)}$ to both sides of this equation gives us (according to (3.3))

$$P_{Z_1} \delta q_\beta^{t_0+0} = \alpha_0(w) \cdot P_{Z_1} v_\beta^{t_0+0}.$$

Hence the lemma follows. □

**Remark 3.7.** In the case $|Z_1| \leq \nu - 2$ we immediately get that the exceptional set (for which $P_{\overline{Z_1}}\left(v_1^{t_0+0} - v_2^{t_0+0}\right) = 0$) has at least two codimensions.

**Lemma 3.8.** *Assume that $\bigcup_{i=1}^n Z_i = \mathcal{A}$ and the advance functionals $\alpha_0, \ldots, \alpha_n$ are equal. Then the orbit segment $\omega$ is sufficient.*

**Proof.** Consider an arbitrary neutral vector $w = \left(\delta q_1^{t_0+0}, \delta q_2^{t_0+0}\right)$ with the advance $\alpha_0(w) = \cdots = \alpha_n(w)$. We need to prove that $w$ is a scalar multiple of the velocity, i. e. $\delta q_\beta^{t_0+0} = \alpha_0(w) \cdot v_\beta^{t_0+0}$ for $\beta = 1, 2$. By subtracting the vector $\alpha_0(w) \cdot \left(v_1^{t_0+0}, v_2^{t_0+0}\right)$ from $w$, we can achieve that $\alpha_0(w) = \cdots = \alpha_n(w) = 0$. Then we need to show that $w = 0$.

The equation $\alpha_0(w) = 0$ implies that $\delta q_1^{t_0+0} = \delta q_2^{t_0+0} := y \in \mathbb{R}^\nu$. Then (3.3) and $\alpha_1(w) = 0$ yield that $R_{Z_1(1)}y = R_{Z_1(2)}y$, i. e. $y = R_{Z_1}y$. This means, however, that $P_{Z_1}y = 0$ and
$$\delta q_1^{t_1-0} = \delta q_2^{t_1-0} = \delta q_1^{t_1+0} = \delta q_2^{t_1+0} = R_{Z_1(1)}y.$$

By continuing this argument we obtain that
$$\delta q_1^{t_i-0} = \delta q_2^{t_i-0} = \delta q_1^{t_i+0} = \delta q_2^{t_i+0} = R_{Z_i(1)} \cdot R_{Z_{i-1}(1)} \cdot \ldots \cdot R_{Z_1(1)} y := y_i$$
for $i = 1, \ldots, n$, and $P_{Z_i}y_{i-1} = 0$. The last equation, however, means that $P_{Z_i}y = 0$ for $i = 1, \ldots, n$, i. e. $P_{\mathcal{A}}y = 0$, because of the assumption $\bigcup_{i=1}^n Z_i = \mathcal{A}$. By the reduction equations we have $P_{\overline{\mathcal{A}}}x = 0$, where $\overline{\mathcal{A}} = \{1, \ldots, \nu\} \setminus \mathcal{A}$. This finishes the proof of the lemma. □

**Lemma 3.9.** *Assume that $|Z_1| = \nu$, $0 < |Z_n| < \nu$, and $Z_2 = Z_3 = \cdots = Z_{n-1} = \emptyset$. We claim that $\alpha_0 = \alpha_1 = \cdots = \alpha_n$ (that is, by virtue of the previous lemma, the orbit segment $\omega$ is sufficient), provided that*

*(i) $P_{\overline{Z_n}}\left(v_1^{t_{n-1}+0} - v_2^{t_{n-1}+0}\right) \neq 0$, and*

*(ii) $P_{Z_n}\left(v_1^{t_0+0} - v_2^{t_0+0}\right) \neq 0$.*

**Proof.** Consider an arbitrary neutral vector $w \in \mathcal{N}(\omega)$ with the configuration variations $\delta q_\beta^{t_i \pm 0}$, $i = 0, \ldots, n$, $\beta = 1, 2$. It follows immediately from Lemma 3.6 and from our assumption (i) that $\alpha_1(w) = \cdots = \alpha_n(w) := \alpha$. We remind the reader that the non-equality $P_{\overline{Z_i}}\left(v_1^{t_{i-1}+0} - v_2^{t_{i-1}+0}\right) \neq 0$ always holds for $i = 2, 3, \ldots, n-1$, for $Z_i = \emptyset$ and the relative velocity $v_1^{t_{i-1}+0} - v_2^{t_{i-1}+0}$ is never zero. Moreover, Lemma 3.6 also claims that $P_{Z_n}\delta q_\beta^{t_{n-1}+0} = \alpha \cdot P_{Z_n} v_\beta^{t_{n-1}+0}$ for $\beta = 1, 2$. The velocity reflection equations for the collisions $\sigma_{n-1}, \ldots, \sigma_1$, the hypothesis $Z_2 = \cdots = Z_{n-1} = \emptyset$, and the facts $\alpha_1(w) = \cdots = \alpha_{n-1}(w) = \alpha$ together imply that $P_{Z_n}\delta q_\beta^{t_1-0} = \alpha \cdot P_{Z_n} v_\beta^{t_1-0}$ for $\beta = 1, 2$. By applying the reflection $R_{Z_1(\beta)}$ to the last equation we get that $P_{Z_n}\delta q_\beta^{t_0+0} = \alpha \cdot P_{Z_n} v_\beta^{t_0+0}$ for $\beta = 1, 2$, especially
$$P_{Z_n}\left(\delta q_1^{t_0+0} - \delta q_2^{t_0+0}\right) = \alpha \cdot P_{Z_n}\left(v_1^{t_0+0} - v_2^{t_0+0}\right).$$
This means, however, that $\alpha_0(w) = \alpha$, as long as the assumption (ii) holds. Hence the lemma follows. □

By putting together lemmas 3.6, 3.8, and 3.9, we immediately obtain the proof of Key Lemma 3.5. □

The last two "transversality" lemmas in this section are the analogues of Lemma 4.1 of [10], and they will be used later in Section 5.

**Definition.** Define the set of axes $Z_0 \subset \mathcal{A}$ quite similarly to the definition of $Z_1, \ldots, Z_n$: It is the set of all axes $j \in \mathcal{A}$ for which the number of the $j$-wall collisions $((q)_j = 0$ or $(q)_j = 1)$ of $q_1$ in the time interval $(a, t_0)$ has different parity than the number of the $j$-wall collisions of $q_2$ during the same time interval.

**Lemma 3.10.** *Denote by* $\mathbf{M}(\Sigma; a, b)$ *the set of all phase points* $x_0 \in \mathbf{M}$ *for which* $\omega = S^{[a,b]}x_0$ *is a non-singular orbit segment with the symbolic collision sequence*

$$\Sigma = (Z_0, \sigma_0, Z_1, \sigma_1, \ldots, Z_n, \sigma_n),$$

*and* $a, b$ *are not moments of collision. Let* $j_1, j_2 \in \{1, \ldots, \nu\}$ *be two indices, and consider the following submanifolds* $J_1$ *and* $J_2$ *of* $\mathbf{M}(\Sigma; a, b)$:

$$J_1 = \left\{ x \in \mathbf{M}(\Sigma; a, b) \,\big|\, (v_1^a)_{j_1} - (v_2^a)_{j_1} = 0 \right\},$$

$$J_2 = \left\{ x \in \mathbf{M}(\Sigma; a, b) \,\big|\, (v_1^b)_{j_2} - (v_2^b)_{j_2} = 0 \right\}.$$

*Assume that either* $j_1 \notin Z_0$ *or* $j_1 \in \bigcup_{l=1}^n Z_l$.
*We claim that* $\operatorname{codim}(J_1 \cap J_2) \geq 2$.

**Proof.** The proof to be presented here is going to be a local, geometric argument, pretty similar to those in the proofs of Lemma 3.25 of [17] or Lemma 4.1 of [10]. Without restricting the generality, we can assume that $a = 0$.

Consider an arbitrary point $x_0 \in J_1$, and select a small number $\epsilon_0 > 0$. Denote by $U_0$ the $\epsilon_0$-neighborhood of the base point $x_0$ in $J_1$:

$$U_0 = \{x \in J_1 \,|\, d(x_0, x) < \epsilon_0\}.$$

Besides $\epsilon_0$, select another, small positive number $\delta_0$ such that for every $x \in U_0$ the trajectory segment $S^{[0,\delta_0]}x$ does not have a collision, not even a non-dispersive ball-to-wall collision.

We will foliate $U_0$ by *convex, local orthogonal manifolds* (see Section 2) as follows: We introduce the equivalence relation $\sim$ in $U_0$ by the formula

$$x \sim y \Leftrightarrow (\exists \lambda \in \mathbb{R}) \quad q_1^0(x) - q_1^0(y) = q_2^0(y) - q_2^0(x) = \lambda \cdot e_{j_1} \qquad (3.11)$$

for $x, y \in U_0$, where $e_{j_1}$ is the standard unit vector in the positive direction of the $j_1$-st coordinate axis. It is easy to see that, indeed, the equivalence classes $C(x) \subset U_0$ of $\sim$ are $(\nu + k - 1)$–dimensional, smooth, connected submanifolds of $U_0$. Furthermore, the positive images $S^t(C(x))$ ($0 < t < \delta_0$) of these submanifolds under the action of the flow $S^t$ are local, convex orthogonal manifolds. Actually, one easily sees that the projection $\pi[S^t(C(x))]$ of $S^t(C(x))$ into the configuration space $\mathbf{Q}$ ($x \in U_0$, $0 < t < \delta_0$) is a cylindric hypersurface in $\mathbf{Q}$ with the one-dimensional generator (constituent) space $\{\lambda \cdot (e_{j_1}, -e_{j_1}) \,|\, \lambda \in \mathbb{R}\}$.

Assume now that the base point $x_0 \in U_0$ also belongs to the manifold $J_2$, i.e. $(v_1^b(x_0))_{j_2} - (v_2^b(x_0))_{j_2} = 0$. In order to prove the lemma it is enough to show that the equivalence class $C(x_0)$ of $x_0$ is transversal to the manifold $J_2$. This will be proved, as long as we can show that the image $S^b(C(x_0))$ of $C(x_0)$ is strictly of the phase points $S^b y$ ($y \in C(x_0)$, $d(x_0, y) \ll 1$) vary with the maximum rank $d - 1 = \nu + k - 1$ in the sphere of velocities $\mathbb{S}^{d-1}$.

In order to prove the strict convexity, however, it is enough to show that the single flat direction $\{\lambda \cdot (e_{j_1}, -e_{j_1}) | \lambda \in \mathbb{R}\}$ of $C(x_0)$ "acquires" a positive curvature during the collisions $\sigma_0, \sigma_1, \ldots, \sigma_n$, i. e. this flat direction does not survive $S^b$ as a neutral direction, see also Section 2. Denote by $l_0$ the smallest number $l \in \{0, 1, \ldots, n\}$ for which the following property holds: Either $l = 0$ and $j_1 \notin Z_0$, or $l > 0$ and $j_1 \in Z_l$. Such an index $l_0$ exists by the assumption of the lemma.

By using an induction on the integers $l = 0, 1, \ldots, l_0 - 1$, easy calculation proves that for every non-negative integer $l < l_0$ the images

$$D\left(S^{t_l - 0}\right)[(e_{j_1}, -e_{j_1}; 0, 0)] \quad \text{and} \quad D\left(S^{t_l + 0}\right)[(e_{j_1}, -e_{j_1}; 0, 0)]$$

of the tangent vector $\delta q_1 = e_{j_1}$, $\delta q_2 = -e_{j_1}$, $\delta v_i = 0$ under the linearization of the flow are equal to the same vector $\pm(e_{j_1}, e_{j_1}; 0, 0)$. Especially, the initial tangent vector $(e_{j_1}, -e_{j_1}; 0, 0)$ remains neutral (with the advance zero) with respect to the collision $\sigma_l$. We claim that this initial tangent vector $(e_{j_1}, -e_{j_1}; 0, 0) \in \mathcal{T}_{x_0} \mathbf{M}$ can not be neutral with respect to the collision $\sigma_{l_0}$. Indeed, the neutrality of the image tangent vector

$$D\left(S^{t_{l_0} - 0}\right)[(e_{j_1}, -e_{j_1}; 0, 0)] = \pm(e_{j_1}, -e_{j_1}; 0, 0)$$

with respect to the collision $\sigma_{l_0}$ means that

$$v_1^{t_{l_0} - 0}(x_0) - v_2^{t_{l_0} - 0}(x_0) = \lambda \cdot e_{j_1} \quad (\lambda \in \mathbb{R}).$$

On the other hand, the vector $(e_{j_1}, -e_{j_1})$ must be orthogonal to the velocity

$$\left(v_1^{t_{l_0} - 0}(x_0), v_2^{t_{l_0} - 0}(x_0)\right).$$

This means, however, that $v_1^{t_{l_0} - 0}(x_0) - v_2^{t_{l_0} - 0}(x_0) = 0$, a contradiction. Therefore, the image vector $D\left(S^{t_{l_0} - 0}\right)[(e_{j_1}, -e_{j_1}; 0, 0)]$ can not be neutral with respect to $\sigma_{l_0}$.

This proves the lemma. □

**Remark 3.12.** It follows easily from the above proof that even if we drop the combinatorial hypotheses on $j_1$ ($j_1 \notin Z_0$ or $j_1 \in \bigcup_{l=1}^n Z_l$), the manifolds $J_1$ and $J_2$ can (locally) coincide only if $j_1 \in Z_0 \cap Z_{n+1}$ and $j_1 \notin \bigcup_{l=1}^n Z_l$. In that case, however, the manifolds $J_1$ and $J_2$ are obviously identical.

**Lemma 3.13.** *Let* $a = 0$, $j \in \{1, \ldots, \nu\}$ *and, similarly to the previous lemma, introduce the manifold*

$$\hat{J}_1 = \left\{x \in \mathbf{M} \mid \left(v_1^0\right)_j - \left(v_2^0\right)_j = 0\right\}.$$

*Assume that* $x_0 \in \hat{J}_1 \cap \mathrm{int}\mathbf{M}$ *is a smooth point of* $\hat{J}_1$ *for which the positive orbit* $S^{(0, \infty)} x_0$ *contains a singularity at* $t^* = t^*(x_0) > 0$, *and* $x_0$ *is also a smooth point of the corresponding singularity set* $\mathcal{S} \ni x_0$. *We claim that* $\hat{J}_1$ *and* $\mathcal{S}$ *intersect each other transversally at* $x_0$.

**Proof.** Following the proof of the previous lemma, we again consider the foliation $U_0 \cap \hat{J}_1 = \bigcup_{\alpha \in A} C(x_\alpha)$, based on the equivalence relation (3.11), where $A$ is a suitable index set. As we have seen in the previous proof, for every small $t > 0$ the set $S^t(C(x_\alpha))$ is, indeed, a local, convex, orthogonal manifold. Sublemma 4.2 of [10] then yields that the manifold $C(x_0) \ni x_0$ is transversal to $\mathcal{S}$ at the base point $x_0$.

Hence the lemma follows. □

## 4. Richness Is Abundant

In this section we will be investigating *non-singular* trajectories

$$S^{(-\infty,\infty)} x_0 = \{x_t = (q_1^t, q_2^t; \dot{q}_1^t, \dot{q}_2^t) \mid t \in \mathbb{R}\}$$

of the standard $(\nu, k, r)$-flow with $k > 0$. (Sometimes we supress the time derivatives $\dot{q}_i^t$, and just simply write $x_t = (q_1^t, q_2^t)$.) The singular orbits will be taken care of in Section 5.

First of all, we note that, according to the theorems by Vaserstein [25] and Galperin [6], the moments of collisions on the trajectory $S^{(-\infty,\infty)} x_0$ can not accumulate at a finite point, i. e. every bounded time interval only contains finitely many collisions.

The next thing we want to make sure is that there are infinitely many dispersive (i. e. "ball-to-ball") collisions both in $S^{[0,\infty)} x_0$ and in $S^{(-\infty,0]} x_0$.

**Lemma 4.1.** *There exist finitely many nonzero vectors $z_1, \ldots, z_l$ in $\mathbb{R}^{\nu-k}$ with the following property:*

*If a non-singular trajectory $S^{(-\infty,\infty)} x_0$ has no ball-to-ball collision on a time interval $(t_0, \infty)$, then either*

*(1) $x_0$ belongs to some slim exceptional subset of the phase space $\mathbf{M}$,*

*or*

*(2) the projection $P_{\overline{A}}\left(v_1^{t_0} - v_2^{t_0}\right) \in \mathbb{R}^{\nu-k}$ of the relative velocity (measured at $t_0$) into the component $\mathbb{R}^{\nu-k}$ of the velocity space $\mathbb{R}^\nu = \mathbb{R}^k \oplus \mathbb{R}^{\nu-k} = \mathbb{R}^{\mathcal{A}} \oplus \mathbb{R}^{\overline{\mathcal{A}}}$ is perpendicular to some vector $z_i$.*

*In either case the entire trajectory $S^{(-\infty,\infty)} x_0$ does not contain a single ball-to-ball collision.*

**Remark.** In the case $k = \nu$ we have that $l = 0$ and the second possibility (2) will not occur.

**Proof.** Assume that $t_0 = 0$, i. e. the positive orbit

$$S^{[0,\infty)} x_0 = \{x_t = (q_1^t, q_2^t) \mid t \geq 0\}$$

does not have any ball-to-ball collision. Then there is a quite standard method of "unfolding" this positive orbit by reflecting the container $\mathbf{C} = [0,1]^k \times \mathbb{T}^{\nu-k}$ across its boundary hyperplanes $(q)_j = 0$, $(q)_j = 1$, $j = 1, \ldots, k$. Namely, we select first an arbitrary Euclidean lifting $\tilde{q}_1^0, \tilde{q}_2^0 \in \mathbb{R}^\nu$ of the initial positions $q_1^0, q_2^0 \in [0,1]^k \times \mathbb{T}^{\nu-k}$ extend this lifting to $t \geq 0$ in a linear manner: $\tilde{q}_i^t := \tilde{q}_i^0 + tv_i^0$, $t \geq 0$, $i = 1, 2$. (We note that $(\tilde{q}_1^t, \tilde{q}_2^t)$ is *not* a lifting of the original orbit.) Finally, we project the linear orbit $\{(\tilde{q}_1^t, \tilde{q}_2^t)|\, t \geq 0\}$ into the torus

$$\left(\mathbb{R}^k/2 \cdot \mathbb{Z}^k\right) \times \mathbb{T}^{\nu-k} = \mathbb{R}^\nu / \left[(2 \cdot \mathbb{Z}^k) \times \mathbb{Z}^{\nu-k}\right]$$

by using the natural projection

$$\mathbb{R}^\nu \longrightarrow \mathbb{R}^\nu / \left[(2 \cdot \mathbb{Z}^k) \times \mathbb{Z}^{\nu-k}\right],$$

and obtain the positive trajectory $\{(\hat{q}_1^t, \hat{q}_2^t)|\, t \geq 0\}$. It is obvious that the finally obtained orbit is independent of the initial selection of the lifting. Let us define the canonical "folding map"

$$\Phi = (\phi, \ldots, \phi;\, \mathrm{id}, \ldots, \mathrm{id}) : \left(\mathbb{R}^k/2 \cdot \mathbb{Z}^k\right) \times \mathbb{T}^{\nu-k} \to [0,1]^k \times \mathbb{T}^{\nu-k} = \mathbf{C}$$

in such a way that the same 2-periodic rooftop function $\phi : \mathbb{R}/2 \cdot \mathbb{Z} \to [0,1]$, $\phi(x) := d(x, 2 \cdot \mathbb{Z})$, acts on the first $k$ components, while the identity function acts on the remaining coordinates. One easily sees that $\Phi\left((\hat{q}_1^t, \hat{q}_2^t)\right) = (q_1^t, q_2^t) = x_t$.

The fact that the original positive orbit $S^{[0,\infty)}x_0$ has no ball-to-ball collision immediately implies these assertions:

(i) $d(\hat{q}_1^t, \hat{q}_2^t) > 2r$ for all $t \geq 0$;

(ii) $d\left(\hat{q}_1^t, R_{\mathcal{A}}\hat{q}_2^t\right) > 2r$ for all $t \geq 0$,

where $d(.,.)$ is the usual Euclidean distance in $\left(\mathbb{R}^k/2 \cdot \mathbb{Z}^k\right) \times \mathbb{T}^{\nu-k}$ inherited from $\mathbb{R}^\nu$ and the operator $R_{\mathcal{A}}$ multiplies the first $k$ coordinates by minus one, see also the previous section. According to Lemma A.2.2 of [23], there are finitely many nonzero vectors $w_1, \ldots, w_p \in \mathbb{R}^\nu$ (none of which is a scalar multiple of another) — not depending on the phase point $x_0$ — such that $\langle v_1^0 - v_2^0;\, w_{j_1}\rangle = 0$ for some $j_1$ (because of (i)), and $\langle v_1^0 - R_{\mathcal{A}} v_2^0;\, w_{j_2}\rangle = 0$ for some $j_2$ (because of (ii)). If at least one of the two projections $P_{\mathcal{A}} w_{j_1}$ and $P_{\mathcal{A}} w_{j_2}$ is nonzero, then, as it follows from the independence of the vectors $P_{\mathcal{A}}(v_1^0 - v_2^0)$ and $P_{\mathcal{A}}(v_1^0 - R_{\mathcal{A}} v_2^0) = P_{\mathcal{A}}(v_1^0 + v_2^0)$, the equations $\langle v_1^0 - v_2^0;\, w_{j_1}\rangle = 0$ and $\langle v_1^0 - R_{\mathcal{A}} v_2^0;\, w_{j_2}\rangle = 0$ are independent, so they together define a manifold with codimension two. Such phase points $x_0$ are listed up in part (1) of the lemma.

Therefore, we can assume that $w_{j_1}$ and $w_{j_2}$ both belong to the second factor $\mathbb{R}^{\nu-k}$ of the velocity space $\mathbb{R}^\nu = \mathbb{R}^k \oplus \mathbb{R}^{\nu-k}$. If $w_{j_1} \neq w_{j_2}$ (i. e. they are not even parallel), then the equations $\langle v_1^0 - v_2^0;\, w_{j_1}\rangle = 0$ and $\langle v_1^0 - R_{\mathcal{A}} v_2^0;\, w_{j_2}\rangle = 0$ are again independent, and such phase points $x_0$ are listed up in (1).

The only way of obtaining a codimension-one family of exceptional phase points $x_0$ is then to have $w_{j_1} = w_{j_2}$ and $P_{\mathcal{A}} w_{j_1} = 0$. List up the vectors $w_j$ with $P_{\mathcal{A}} w_j = 0$ in a sequence $z_1, \ldots, z_l$ ($\in \mathbb{R}^{\nu-k}$), and obtain the first statement of the lemma.

The fact that for an exceptional phase point $x_0$ – listed up in (1)-(2) – the entire trajectory $S^{(-\infty,\infty)}x_0$ does not contain a single ball-to-ball collision immediately follows from the following, simple observation:

**Sublemma 4.2.** *Consider two opposite rays $L^+$, $L^- \subset \mathbb{R}^n$ as follows:*

$$L^+ = \{q_0 + tv_0|\, t \geq 0\}, \quad L^- = \{q_0 + tv_0|\, t \leq 0\},$$

$q_0 \in \mathbb{R}^n$, $0 \neq v_0 \in \mathbb{R}^n$. *Then $d(L^+, \mathbb{Z}^n) = d(L^-, \mathbb{Z}^n)$.*

**Proof.** Apply the natural projection $\pi : \mathbb{R}^n \to \mathbb{R}^n/\mathbb{Z}^n = \mathbb{T}^n$ to the half lines $L^\pm$. It is well known that $\text{Cl}(\pi L^+) = \text{Cl}(\pi L^-) = S$ is a coset with respect to a subtorus of $\mathbb{T}^n$. (Here Cl denotes the closure.) This means that if $d(L^+, \mathbb{Z}^n) < \alpha$, then $S$ intersects the open ball with radius $\alpha$ centered at zero in $\mathbb{T}^n$ and, therefore, $d(L^-, \mathbb{Z}^n) < \alpha$, as well.

This finishes the proof of Sublemma 4.2 and, hence, the proof of Lemma 4.1, as well. □

**Remark.** Part (2) of the assertion of 4.1 is really not vague. Indeed, in the case of a trajectory without a ball-to-ball collision the time evolutions of the $(q)_j$-components ($1 \leq j \leq \nu$) are *independent* of each other, thanks to the orthogonality of the walls of the container. Therefore, if the $\overline{\mathcal{A}}$-component $P_{\overline{\mathcal{A}}}(v_1^0 - v_2^0)$ of the initial relative velocity belongs to some exceptional hyperplane of $\mathbb{R}^{\nu-k}$ (just as in the proof of Lemma A.2.2 of [23]) and the the initial relative position $P_{\overline{\mathcal{A}}}(q_1^0 - q_2^0)$ belongs to some non-empty, open region, then even the $\overline{\mathcal{A}}$-parts of the positions $q_i^t$ ($i = 1, 2$) never get closer to each other than $2r$, and there will be no ball-to-ball collision.

Thus, by dropping the exceptional phase points $x_0$ listed up in (1) and (2) of 4.1 (the latter set of phase points will be discussed in the next section), we can assume that the non-singular trajectory $S^{(-\infty,\infty)}x_0$ contains infinitely many ball-to-ball collisions in each time direction. Let us list all such collisions of $S^{(-\infty,\infty)}x_0$ as $(\ldots, \sigma_{-1}, \sigma_0, \sigma_1, \ldots)$, so that $\sigma_0$ is the first collision occuring in positive time ($t = 0$ is supposed not to be a moment of collision), i. e. there is not even a non-dispersive, wall-to-ball collision in the time interval $[0, t_0)$. Just as in the previous section, $t_n = t(\sigma_n)$ denotes the time of the $n$-th collision, $n \in \mathbb{Z}$. Following (3.2), we can now speak about the sets $Z_i(1)$, $Z_i(2)$, $Z_i \subset \mathcal{A}$, ($i \in \mathbb{Z}$). The doubly infinite sequence

$$\Sigma = (\ldots, \sigma_{-1}, Z_0, \sigma_0, Z_1, \sigma_1, \ldots)$$

is now called the symbolic collision sequence of the trajectory $\omega = S^{(-\infty,\infty)}x_0$.

**From now on, in this section we will always be assuming the following inductive hypothesis:**

**Hypothesis 4.3.** *For every integer $k' < k$ the standard $(\nu, k', r)$-flow* **is ergodic**, *and it enjoys the following, additional properties:*

*(A) There exists an exceptional, slim subset $E \subset \mathbf{M}(\nu, k', r)$ of the phase space of this flow such that the trajectory $S^{(-\infty,\infty)}y$ of every phase point $y \in \mathbf{M}(\nu, k', r) \setminus E$ either has no ball-to-ball collisions (such phase points will be taken care of in the next section), or $S^{(-\infty,\infty)}y$ has at most one singular collision and it is also sufficient. (In the singular case this means that both branches are sufficient, see Section 2.)*

*(B) The Chernov–Sinai Ansatz holds for the $(\nu, k', r)$-flow, see Condition 3.1 in [10].*

**Remark.** The notions and basic properties of sufficiency and slimness are summarized in Section 2.

The main result of this section is the analogue of Lemma 3 of [9] or Main Lemma 5.1 of [11]:

**Key Lemma 4.4 (Strong Ball Avoiding Theorem).** *Assume 4.3. Then there exists a slim subset $S$ of the phase space $\mathbf{M} = \mathbf{M}(\nu, k, r)$ with the following property:*

*If the symbolic sequence $\Sigma$ of a non-singular orbit $\omega = S^{(-\infty,\infty)}x_0$ (with infinitely many ball-to-ball collisions) is not rich in the sense of 3.4 (i. e. $\bigcup_{i=-\infty}^{\infty} Z_i \neq \mathcal{A}$ or $|Z_i|(\nu - |Z_i|) = 0$ for every $i \in \mathbb{Z}$), then $x_0 \in S$.*

All the rest of this section will be devoted to the proof of the key lemma. The proof will be split into a few lemmas. The first one of them takes care of the case $\bigcup_{i=-\infty}^{\infty} Z_i \neq \mathcal{A}$.

**Lemma 4.5.** *There exists a slim set $S \subset \mathbf{M}$ such that if $\bigcup_{i=-\infty}^{\infty} Z_i \neq \mathcal{A}$, then the phase point $x_0$ necessarily belongs to $S$.*

**Proof.** We can assume that $1 \notin \bigcup_{i=-\infty}^{\infty} Z_i$. By the basic properties of slim sets (see 2.4–2.12 in this article), it is enough to show that the given phase point $x_0 \in \mathbf{M}$ has an open neighborhood $U_0 \ni x_0$ such that the set $U_0 \cap P_1$ is slim, where

$$P_1 = \left\{ x \in \mathbf{M}^* \,\Big|\, 1 \notin \bigcup_{i=-\infty}^{\infty} Z_i \right\},$$

and $\mathbf{M}^*$ is the set of all non-singular phase points $x \in \mathbf{M}$ for which $S^{(-\infty,\infty)}x$ contains (infinitely many) ball-to-ball collisions. The proof will borrow the main ideas from the proofs of Main Lemma 5.1 and Lemma 5.3 of [11].

Select first a small, open ball neighborhood $U_0 \ni x_0$ of the fixed phase point $x_0 \in P_1$ so that $U_0 \cap \partial \mathbf{M} = \emptyset$. It is an important observation that the orbit $S^{(-\infty,\infty)}x = \{(q_1^t(x), q_2^t(x)) \,|\, t \in \mathbb{R}\}$ of every phase point $x \in U_0 \cap P_1$ can be "unfolded" to obtain another trajectory $\hat{\omega} = \{(\hat{q}_1^t(x), \hat{q}_2^t(x)) \,|\, t \in \mathbb{R}\}$, $\hat{q}_i^t(x) \in (\mathbb{R}/2 \cdot \mathbb{Z}) \times [0,1]^{k-1} \times \mathbb{T}^{\nu-k}$ by again reflecting the original container space $\mathbf{C} = [0,1]^k \times \mathbb{T}^{\nu-k}$ across two of its boundary hyperplanes $(q)_1 = 0$ and $(q)_1 = 1$, much the same way as it was done in the proof of Lemma 4.1. Namely, we again consider the covering map

$$\Phi = (\phi; \mathrm{id}, \ldots, \mathrm{id}) : (\mathbb{R}/2 \cdot \mathbb{Z}) \times [0,1]^{k-1} \times \mathbb{T}^{\nu-k} \to [0,1]^k \times \mathbb{T}^{\nu-k} \qquad (4.6)$$

with the rooftop function $\phi : \mathbb{R}/2 \cdot \mathbb{Z} \to [0,1]$ ($\phi(x) = d(x, 2 \cdot \mathbb{Z})$) acting on the first component, and then pull back the given trajectory

$$S^{(-\infty,\infty)}x = \{(q_1^t(x), q_2^t(x)) \,|\, t \in \mathbb{R}\}$$

by $\Phi$ to obtain the orbit $\hat{\omega} = \{(\hat{q}_1^t(x), \hat{q}_2^t(x)) \,|\, t \in \mathbb{R}\}$, $\hat{q}_i^t(x) \in (\mathbb{R}/2 \cdot \mathbb{Z}) \times [0,1]^{k-1} \times \mathbb{T}^{\nu-k}$, after having fixed a continuous initial pull-back

$$(\hat{q}_1^0(y), \hat{q}_2^0(y)) := (q_1^0(y), q_2^0(y))$$

(mod $2 \cdot \mathbb{Z}$ in the first coordinates) for all phase points $y \in U_0$. It is an important consequence of the assumption $1 \notin \bigcup_{i=-\infty}^{\infty} Z_i(x)$ ($x \in U_0 \cap P_1$) that the curve $\hat{\omega}$ is a two-ball billiard trajectory in the container $\hat{\mathbf{C}} = (\mathbb{R}/2 \cdot \mathbb{Z}) \times [0,1]^{k-1} \times \mathbb{T}^{\nu-k}$, i. e. there are only collisions with the cylinder $d(\hat{q}_1, \hat{q}_2) = 2r$ and collisions never occur at the "antipodal" cylinder $d(\hat{q}_1, R_{\{1\}}\hat{q}_2) = 2r$. In other words, $d(\hat{q}_1^t(x), R_{\{1\}}\hat{q}_2^t(x)) > 2r$ and $d(\hat{q}_1^t(x), \hat{q}_2^t(x)) \geq 2r$ for every $t \in \mathbb{R}$, while $d\left(\hat{q}_1^{t_n(x)}(x), \hat{q}_2^{t_n(x)}(x)\right) = 2r$ ($n \in \mathbb{Z}$), where $t_n(x) = t(\sigma_n(x))$ is the time of the $n$-th dispersive collision $\sigma_n$ on the orbit of $x \in U_0 \cap P_1$.

The above discussed $\hat{\mathbf{C}}$-dynamics $\{(\hat{q}_1^t(x), \hat{q}_2^t(x)) \mid t \in \mathbb{R}\}$ (which is a standard $(\nu, k-1, r)$-flow without the normalizations $(v_1)_1 + (v_2)_1 = (q_1)_1 + (q_2)_1 = 0$) can now be defined for every phase point $x \in U_0$ by using the initial lifting $(\hat{q}_1^0(x), \hat{q}_2^0(x)) = (q_1^0(x), q_2^0(x))$ of the positions, irrespectively of whether $x \in P_1$ or not. (In other words, when defining this $\hat{\mathbf{C}}$-dynamics $\{(\hat{q}_1^t(x), \hat{q}_2^t(x)) \mid t \in \mathbb{R}\}$, we do not remove the antipodal cylinder $d(\hat{q}_1, R_{\{1\}}\hat{q}_2) < 2r$ from the configuration space but, rather, we allow the above inequalities.) We choose a small number $\epsilon_0 > 0$ and define the following closed subsets of $U_0$ (see also (5.2) and (5.7) of [11]):

$$F_+ = \{x \in U_0 \mid d(\hat{q}_1^t(x), R_{\{1\}}\hat{q}_2^t(x)) \geq 2r \quad \forall t \geq 0\},$$
$$F_- = \{x \in U_0 \mid d(\hat{q}_1^t(x), R_{\{1\}}\hat{q}_2^t(x)) \geq 2r \quad \forall t \leq 0\},$$
(4.7)

$$F_+' = \{x \in U_0 \mid d(\hat{q}_1^t(x), R_{\{1\}}\hat{q}_2^t(x)) \geq 2r - \epsilon_0 \quad \forall t \geq 0\},$$
$$F_-' = \{x \in U_0 \mid d(\hat{q}_1^t(x), R_{\{1\}}\hat{q}_2^t(x)) \geq 2r - \epsilon_0 \quad \forall t \leq 0\},$$
(4.8)

where, in the case of singular trajectories, we understand these inequalities in such a way that they should hold for *some* trajectory branch. This convention makes the sets $F_\pm \subset U_0$, $F_\pm' \subset U_0$ closed, just as in Main Lemma 5.1 of [11] or in (5.4)–(5.7) of [12]. It is obvious that

$$F_+ \subset F_+', \qquad F_- \subset F_-',$$
(4.9)

and

$$U_0 \cap P_1 \subset F_- \cap F_+.$$
(4.10)

Since the $\hat{\mathbf{C}}$-flow $\{(\hat{q}_1^t(x), \hat{q}_2^t(x)) \mid t \in \mathbb{R}\}$ has the quantities $(v_1)_1 + (v_2)_1$ and $(q_1)_1 + (q_2)_1$ as first integrals, it is now quite natural to define a foliation

$$U_0 = \cup \{H_{a,b} \mid a \in I_1, \ b \in I_2\}$$
$$H_{a,b} := \{x \in U_0 \mid (q_1(x))_1 + (q_2(x))_1 = a, \ (v_1(x))_1 + (v_2(x))_1 = b\},$$
(4.11)

where $I_1, I_2 \subset \mathbb{R}$ are suitable open intervals. (We can assume that the shape of the small, open neighborhood $U_0 \ni x_0$ is such that it permits us to establish the union (4.11) with open intervals $I_1, I_2$.)

It is an important consequence of the assumption 4.3/(A) that for $\mu$-almost every phase point $x \in U_0$ there exist the maximum-dimensional (actually, $(\nu + k - 2)$-dimensional) local, exponentially stable and unstable invariant manifolds $\gamma^s(x)$, $\gamma^u(x)$ (containing $x$ as an interior point) with respect to the $\hat{\mathbf{C}}$-dynamics, defineable as follows:

$$\gamma^s(x) = \mathrm{CC}_x \{y \in U_0 \mid d(\hat{x}^t, \hat{y}^t) \to 0 \text{ exp. fast as } t \to +\infty\},$$
$$\gamma^u(x) = \mathrm{CC}_x \{y \in U_0 \mid d(\hat{x}^t, \hat{y}^t) \to 0 \text{ exp. fast as } t \to -\infty\},$$
(4.12)

where

$$\hat{x}^t = (\hat{q}_1^t(x), \hat{q}_2^t(x); \hat{v}_1^t(x), \hat{v}_2^t(x)),$$

and the $\hat{\mathbf{C}}$-phase point $\hat{y}^t$ is defined analogously. (Here the symbol $\mathrm{CC}_x(.)$ denotes the connected component of a set containing the point $x$.) Clearly, $\gamma^s(x) \cup \gamma^u(x) \subset$

According to the fundamental "integrability" of slimness (see Section 2 or Lemma 3.8 of [16]), it is enough to show that the closed subset $F_- \cap F_+ \cap H_{a,b}$ of $H_{a,b}$ is slim in $H_{a,b}$ for every $a \in I_1$, $b \in I_2$. It follows immediately from the ergodicity assumption of 4.3/A that

$$\mu_{H_{a,b}}\left((F'_- \cup F'_+) \cap H_{a,b}\right) = 0. \tag{4.13}$$

Recall that $F'_\pm \supset F_\pm$. In order to show that the closed, zero set $F_- \cap F_+ \cap H_{a,b}$ has at least two codimensions in $H_{a,b}$, it is enough to prove that this set enjoys the non-separating property, see (ii) of Lemma 2.9 here, or (ii) of Lemma 3.9 in [16]. This is what we will do.

Besides the exponentially stable and unstable manifolds $\gamma^s(x)$, $\gamma^u(x) \subset H_{a,b}$ ($x \in H_{a,b}$), we need to use the one-dimensional neutral manifolds $\gamma^0(x) \subset H_{a,b}$ corresponding to the $\hat{\mathbf{C}}$-flow direction. We set

$$\gamma^0(x) = \Big\{y \in H_{a,b} \,\Big|\, v_i(y) = v_i(x) \text{ and } \exists \lambda \in \mathbb{R} \text{ such that}$$
$$q_i(y) = q_i(x) + \lambda v_i(x) - \frac{\lambda b}{2}(1, 0, \ldots, 0) \quad (i = 1, 2)\Big\} \tag{4.14}$$

for $x \in H_{a,b}$, $a \in I_1$, $b \in I_2$. Subtracting the vector $\dfrac{\lambda b}{2}(1, 0, \ldots, 0)$ from the positions is just needed in order to project back into $H_{a,b}$. It is clear that the following transversality condition holds, see also (5.19) and 5.20 of [12]:

$$\mathcal{T}_x H_{a,b} = \mathcal{T}_x \gamma^u(x) + \mathcal{T}_x \gamma^s(x) + \mathcal{T}_x \gamma^0(x), \tag{4.15}$$

where $\mathcal{T}_x(.)$ denotes the tangent space of a manifold at the foot point $x$ and $+$ is a notation for the (not necessarily orthogonal) direct sum of linear spaces.

It is easy to see that for every $x \in U_0$, $y \in \gamma^s(x)$ or $y \in \gamma^0(x)$

$$d\left(\hat{q}^t(x), \hat{q}^t(y)\right) \leq d\left(\hat{q}^0(x), \hat{q}^0(y)\right) \leq \operatorname{diam} U_0 \text{ for } t \geq 0, \tag{4.16}$$

and, analogously, for any pair $x, y \in U_0$ for which $y \in \gamma^u(x)$ or $y \in \gamma^0(x)$

$$d\left(\hat{q}^t(x), \hat{q}^t(y)\right) \leq d\left(\hat{q}^0(x), \hat{q}^0(y)\right) \leq \operatorname{diam} U_0 \text{ for } t \leq 0. \tag{4.17}$$

Therefore, if the size diam$U_0$ is chosen small enough compared to $\epsilon_0$, then we have the analogue of Lemma 5.8 from [11]:

**Sublemma 4.18.** *For any $x \in U_0$, if $\gamma^s(x) \cap F_+ \neq \emptyset$ ($\gamma^0(x) \cap F_+ \neq \emptyset$), then $\gamma^s(x) \subset F'_+$ ($\gamma^0(x) \subset F'_+$). Analogously, if $\gamma^u(x) \cap F_- \neq \emptyset$ ($\gamma^0(x) \cap F_- \neq \emptyset$), then $\gamma^u(x) \subset F'_-$ ($\gamma^0(x) \subset F'_-$).*

According to our assumption 4.3, if the fixed base point $x_0$ does not belong to some exceptional, slim subset of $\mathbf{M}$ determined by the $\hat{\mathbf{C}}$-dynamics, then $x_0$ is sufficient with respect to this $\hat{\mathbf{C}}$-flow and the fundamental statement of Lemma 3 of [22] (or Theorem 3.6 of [10]) holds true. (Since their formalism is too technical, we do not even quote them here.) As a direct corollary of these results, the absolute continuity of the triple of foliations $\gamma^s(.)$, $\gamma^u(.)$, $\gamma^0(.)$ (see Theorem 4.1 of [7]), and the transversality relation (4.15), we obtain the crucial Zig-zag Lemma, the precise analogue of Corollary 3.10 of [10]:

**Sublemma 4.19 (Zig-zag Lemma).** *Suppose that the base point $x_0 \in \mathbf{M}$ is sufficient with respect to the $\hat{\mathbf{C}}$-dynamics. Then we can select the open neighborhood $U_0 \ni x_0$ so small that the following assertion holds true:*

*For every manifold $H_{a,b} \subset U_0$ ($a \in I_1$, $b \in I_2$), for every open, connected subset $G \subset H_{a,b}$, and for almost every pair of points $x, y \in G \setminus (F'_- \cup F'_+)$ one can find a finite sequence*

$$\xi_1^u,\ \xi_1^0,\ \xi_1^s,\ \xi_2^u,\ \xi_2^0,\ \xi_2^s, \ldots, \xi_n^u,\ \xi_n^0,\ \xi_n^s$$

*of submanifolds of $G$ with the following properties:*

*(i) $\xi_i^\alpha$ is an open, connected subset of $\gamma^\alpha(z)$ for some $z \in \xi_i^\alpha$, $\alpha = u, 0, s$, $i = 1, \ldots, n$;*

*(ii) $x \in \xi_1^u$, $y \in \xi_n^s$;*

*(iii)*

$$\emptyset \neq \xi_i^u \cap \xi_i^0 \subset G \setminus (F'_- \cup F'_+),\ \emptyset \neq \xi_i^0 \cap \xi_i^s \subset G \setminus (F'_- \cup F'_+),\ \emptyset \neq \xi_i^s \cap \xi_{i+1}^u \subset G \setminus (F'_- \cup F'_+)$$

*for $i = 1, \ldots, n$.*

**Remark.** The non-empty intersections in (iii) must be one-point-sets, according to the transversality (4.15).

It follows immediately from (i)–(iii) and Sublemma 4.18 that

$$\bigcup \{\xi_i^\alpha \mid \alpha = u, 0, s;\quad i = 1, \ldots, n\} \subset G \setminus (F_- \cap F_+). \tag{4.20}$$

Thus we obtained that almost every pair of points $x, y \in G \setminus (F'_- \cup F'_+)$ can be connected by a continuous curve $\{\xi(t) \mid 0 \leq t \leq 1\}$ so that $\xi(t) \in G \setminus (F_- \cap F_+)$ for $0 \leq t \leq 1$. Since the open subset $G \setminus (F'_- \cup F'_+)$ has full measure in $G$ (and, therefore, it is dense in $G$), we get that the open subset $G \setminus (F_- \cap F_+)$ of $G$ must be connected. This proves (ii) of Lemma 3.9 of [16], i. e. $\dim(F_- \cap F_+ \cap H_{a,b}) \leq \dim H_{a,b} - 2$ and, therefore, $F_- \cap F_+$ is indeed a slim subset of $U_0$.

The proof of Lemma 4.5 is now complete. $\square$

The last lemma of this section takes care of the case $k = \nu$, $|Z_i|\,(\nu - |Z_i|) = 0$, ($i \in \mathbb{Z}$).

**Lemma 4.21.** *Assume that $k = \nu$ and $|Z_i|\,(\nu - |Z_i|) = 0$ for every integer $i$. Then the phase point $x_0$ necessarily belongs to some slim, exceptional subset of $\mathbf{M}$ not depending on $x_0$.*

**Proof (A sketch).** Since the proof of this lemma is very similar to that of Lemma 4.5, we are not going to present it here in whole detail but, instead, we will just point out the differences between the two approaches.

We again begin with "unfolding" the trajectory

$$S^{(-\infty,\infty)} x_0 = \{(q_1^t(x_0), q_2^t(x_0)) \mid t \in \mathbb{R}\} \tag{4.22}$$

by reflecting the cubic container $\mathbf{C} = [0,1]^\nu$ across its faces $(q)_j = 0$ and $(q)_j = 1$, $j = 1, \ldots, \nu$, pretty much the same way as we did in the proof of 4.5. Namely, we consider the covering map with the rooftop map $\phi: \mathbb{R}/2 \cdot \mathbb{Z} \to [0, 1]$, $\phi(y) := d(y, 2 \cdot \mathbb{Z})$ acting on all components, see also (4.6). Then we just pull back the trajectory of (4.22) by the mapping $\Phi$ to obtain the unfolded trajectory

$$\hat{\omega}(x_0) = \left\{ \left(\hat{q}_1^t(x_0), \hat{q}_2^t(x_0)\right) \mid t \in \mathbb{R} \right\},$$

$\hat{q}_i^t(x_0) \in \mathbb{R}^\nu/2 \cdot \mathbb{Z}^\nu$ by using the selected initial pull-back

$$\left(\hat{q}_1^0(x_0), \hat{q}_2^0(x_0)\right) = \left(q_1^0(x_0), q_2^0(x_0)\right), \mod 2 \cdot \mathbb{Z}^\nu.$$

The hypothesis $|Z_i|\,(\nu - |Z_i|) = 0$ $(i \in \mathbb{Z})$ implies that the obtained unfolded curve

$$\hat{\omega} \subset (\mathbb{R}^\nu/2 \cdot \mathbb{Z}^\nu) \times (\mathbb{R}^\nu/2 \cdot \mathbb{Z}^\nu) = \hat{\mathbf{C}} \times \hat{\mathbf{C}}$$

has reflections only at the cylinder

$$C = \{(\hat{q}_1, \hat{q}_2) \mid d(\hat{q}_1, \hat{q}_2) = 2r\}$$

and the "antipodal" cylinder

$$\overline{C} = \{(\hat{q}_1, \hat{q}_2) \mid d(\hat{q}_1, -\hat{q}_2) = 2r\}.$$

We can now define the unfolded orbit

$$\hat{\omega}(x) = \left\{ \left(\hat{q}_1^t(x), \hat{q}_2^t(x)\right) \mid t \in \mathbb{R} \right\}$$

for every phase point $x = \left(q_1^0(x), q_2^0(x)\right) = \left(\hat{q}_1^0(x), \hat{q}_2^0(x)\right) \in U_0$ by forgetting about the collisions at $d\left(\hat{q}_1, R_Z \hat{q}_2\right) = 2r$ with $0 < |Z| < \nu$ (making these cylinders transparent), and retaining only the collisions at the cylinders $C$ and $\overline{C}$. In the Appendix we show that the so obtained $\hat{\mathbf{C}}$-dynamics ($\hat{\mathbf{C}} = \mathbb{R}^\nu/2 \cdot \mathbb{Z}^\nu$) with the configuration space $\hat{\mathbf{C}} \times \hat{\mathbf{C}} \setminus (C \cup \overline{C})$ is, in fact, a splitting orthogonal cylindric billiard in the sense of the article [23]. There is now one additional first integral, namely the energy $\|\hat{v}_1 + \hat{v}_2\|^2$. Therefore, the foliation of the small open neighborhood $U_0 \ni x_0$ (the analogue of (4.11)) is now going on by specifying the value of $\|\hat{v}_1^0 + \hat{v}_2^0\|^2 = \|v_1^0 + v_2^0\|^2$:

$$U_0 = \bigcup_{a \in I} H_a, \tag{4.24}$$

$$H_a := \left\{ x \in U_0 \mid \|v_1^0(x) + v_2^0(x)\|^2 = a \right\}.$$

It is shown in the Appendix that, for any fixed value of $\|\hat{v}_1 + \hat{v}_2\|^2$, the $\hat{\mathbf{C}}$-flow is the Cartesian product of two dispersive Sinai–billiards and, therefore, it is Bernoulli, and all but two of its Lyapunov exponents are nonzero. Thus, in this situation we have $\dim \gamma^s(x) = \dim \gamma^u(x) = 2\nu - 2$, $\dim \gamma^0(x) = 2$, $\dim H_a = 4\nu - 2$, and for every $x \in U_0$ the linear direct sum of the tangent spaces of the transversal $\gamma^s(x)$, $\gamma^u(x)$ and $\gamma^0(x)$ is exactly the tangent space of the folium $H_a \ni x$ through $x$, see also (4.15). If the original orbit $\omega(x) = S^{(-\infty,\infty)} x$ of a phase point $x \in U_0$ has the property $|Z_i(x)|\,(\nu - |Z_i(x)|) = 0$ $(i \in \mathbb{Z})$, then this fact is reflected by the $\hat{\mathbf{C}}$-orbit $\hat{\omega}(x)$ in such a way that the configuration point $(\hat{q}_1^t(x), \hat{q}_2^t(x))$ never enters a "forbidden" open region, namely the interiors of the cylinders $d\left(\hat{q}_1, R_Z \hat{q}_2\right) \leq 2r$, same way as in the proof of 4.5. Also, by shrinking the forbidden region mentioned above, we can define the larger closed sets $F'_+$, $F'_-$. Thanks to the ergodicity of the $\hat{\mathbf{C}}$-flow, we again have $\mu_{H_a}\left((F'_- \cup F'_+) \cap H_a\right) = 0$. Plainly, the analogues of sublemmas 4.18–4.19 remain valid in the recent situation. By putting together these ingredients, we see that the closed set $F_- \cap F_+ \cap H_a$ has, indeed, at least two codimensions in $H_a$. Finally, by using again the integrability of slimness (Lemma 3.8 of [16]), we obtain that the points $x \in U_0$ with the property $|Z_i(x)|\,(\nu - |Z_i(x)|) = 0$ ($i \in \mathbb{Z}$) form a slim subset of $U_0$. This finishes the proof of Lemma 4.21 and, therefore, the proof of Key Lemma 4.4, as well. $\square$

## 5. Proof of the Theorem

By using the results of the preceding sections, here we prove our

**Theorem.** *The standard $(\nu, k, r)$-flow ($\nu \geq 2$, $0 \leq k \leq \nu$, $0 < r < 1/4$) is ergodic and, therefore, it is actually a Bernoulli flow, see [3] or [15].*

The proof will be an induction on the number $k = 0, \ldots, \nu$. Since the standard $(\nu, 0, r)$-flow is a classical, dispersive Sinai–billiard, it is known to be ergodic, cf. [21], [2], and [22].

Thus, let us assume that $k \geq 1$ and the theorem has been proved for all smaller values of $k$, i. e. assume Hypothesis 4.3. Besides the key lemmas, the inductive step will use three basic ingredients:

(1) (A) the Theorem on Local Ergodicity by Chernov and Sinai, i. e. Theorem 5 of [22], see also Corollary 3.12 of [10];
(2) (B) the method of proving Theorem 6.1 of [16] and that theorem itself;
(3) (C) the method of connecting the several open domains $\Omega_1, \ldots, \Omega_r$ ($r < \infty$) of $\mathbf{M}$ into which $\mathbf{M} = \mathbf{M}(\nu, k, r)$ is split by the union $\mathcal{F}$ of the codimension-one, exceptional submanifolds from part (2) of Lemma 4.1; (That method first appeared at the end of the proof of Lemma 4 in [9].)
(4) (D) the key lemmas.

We begin with the application of (B) by proving the following analogue of Theorem 6.1 of [16]:

**Lemma 5.1 (A local result).** *Suppose that Hypothesis 4.3 holds true. Use all conditions and notations of Lemma 3.10 and assume, for the sake of simplifying the notations, that $a = 0$. We assume, additionally, that $j_2 \in \mathcal{A}$, i. e. $j_2 \leq k$. Denote by $\hat{J}_1 \subset J_1$ a (small) smooth, open and connected subset of $J_1$. We define the following, closed subsets of $\hat{J}_1$:*

$$F^{(1)}_+ = \left\{ x \in \hat{J}_1 \,\bigg|\, j_2 \notin \bigcup_{l=n+1}^{\infty} Z_l(x) \right\}, \tag{5.2}$$

$$F^{(2)}_+ = \left\{ x \in \hat{J}_1 \,\big|\, |Z_l(x)|\,(\nu - |Z_l(x)|) = 0 \text{ for all } l > n \right\}$$

(If the positive orbit $S^{[0,\infty)}x$ is singular, then we understand the above requirements in such a way that they should hold for **some** branch of $S^{[0,\infty)}x$, thus making the sets $F_+^{(1)}$, $F_+^{(2)}$ closed in $\hat{J}_1$.)

We claim that both sets $F_+^{(1)}$, $F_+^{(2)}$ have an empty interior in $\hat{J}_1$.

**Remark 5.3.** In the case $k < \nu$ the set $F_+^{(2)}$ is contained in $F_+^{(1)}$ (meaning $\cup_{l>n} Z_l(x) = \emptyset$ for every $x \in F_+^{(2)}$), and the statement about $F_+^{(2)}$ follows from the claim about $F_+^{(1)}$. The result $\text{int}_{\hat{J}_1} F_+^{(2)} = \emptyset$ is only significant in the case $k = \nu$.

**Proof.** Since the proof is very similar to the proof of Lemma 6.1 of [16], here we will only present a rough sketch of it. Furthermore, the proof of $\text{int}_{\hat{J}_1} F_+^{(2)} = \emptyset$ (in the case $k = \nu$) is quite analogous to the proof of $\text{int}_{\hat{J}_1} F_+^{(1)} = \emptyset$ and, therefore, we will only be dealing with the result $\text{int}_{\hat{J}_1} F_+^{(1)} = \emptyset$. (This analogy is similar to the analogy between the proofs of lemmas 4.5 and 4.21.)

Without restricting generality, we assume that $j_2 = 1$. We argue by contradiction, so we suppose that $\text{int}_{\hat{J}_1} F_+^{(1)} \neq \emptyset$ and, by choosing a smaller open piece of the manifold $\hat{J}_1$, we can assume that $F_+^{(1)} = \hat{J}_1$, i. e. the positive orbit $S^{[0,\infty)}x$ of every point $x \in \hat{J}_1$ follows a pretty irregular pattern: $1 \notin \cup_{l=n+1}^{\infty} Z_l(x)$. Select now a base point $x_0 \in \hat{J}_1$ and a small, open ball neighborhood $U_0 \subset \mathbf{M}(\nu, k, r)$ of $x_0$. We can assume that the time point $b > 0$ is chosen in such a way that there is no collision (not even a wall collision) between $t_n(x) = t(\sigma_n(x))$ and $b$ for every $x \in U_0$.

Let us push forward the open set $U_0$ and $U_0 \cap \hat{J}_1$ by the $b$-th power of the $(\nu, k, r)$-flow: $S^b(U_0) = W_0$, $S^b\left(U_0 \cap \hat{J}_1\right) = K_0$, $S^b(x_0) = y_0$. By our assumption, all trajectories $S^{(0,\infty)}y$ follow the irregular pattern

$$1 \notin \bigcup_{l=0}^{\infty} Z_l(y), \quad y \in K_0. \tag{5.4}$$

Therefore, for every $y \in K_0$ we can again" unfold" the positive orbit

$$S^{(0,\infty)}y = \{(q_1^t(y), q_2^t(y)) \mid t \geq 0\}$$

by reflecting our container $\mathbf{C} = [0,1]^k \times \mathbb{T}^{\nu-k}$ across two of its faces $(q)_1 = 0$, $(q)_1 = 1$, and obtain the positive orbit

$$\hat{\omega}(y) = \{\hat{y}^t = (\hat{q}_1^t(y), \hat{q}_2^t(y)) \mid t \geq 0\} \quad (y \in K_0)$$

in the container

$$\hat{\mathbf{C}} = (\mathbb{R}/2 \cdot \mathbb{Z}) \times [0,1]^{k-1} \times \mathbb{T}^{\nu-k},$$

$\hat{q}_i^t(y) \in \hat{\mathbf{C}}$, see also the proof of Lemma 4.5.

Quite similarly to the proof of Lemma 4.5, the above defined $\hat{\mathbf{C}}$-dynamics $\hat{\omega}(y)$, which is essentially a standard $(\nu, k-1, r)$-flow without the normalization $(v_1)_1 + (v_2)_1 = (q_1)_1 + (q_2)_1 = 0$, can now be defined for every phase point $y \in W_0$ by using the initial lifting $(\hat{q}_1^0(y), \hat{q}_2^0(y)) = (q_1^0(y), q_2^0(y))$ of the positions, irrespectively of

Consider now the codimension-two submanifold $\tilde{W}_0 \subset W_0$

$$\begin{aligned}\tilde{W}_0 = \{y \in W_0 | \ & (q_1(y))_1 + (q_2(y))_1 = (q_1(y_0))_1 + (q_2(y_0))_1 \\ \text{and } & (v_1(y))_1 + (v_2(y))_1 = (v_1(y_0))_1 + (v_2(y_0))_1 \},\end{aligned} \tag{5.5}$$

and the intersection $\tilde{K}_0 = K_0 \cap \tilde{W}_0$, corresponding to fixing the values of $(q_1)_1 + (q_2)_1$ and $(v_1)_1 + (v_2)_1$ according to what these values are for the base point $y_0 = S^b x_0 \in K_0$. The phase point $y \in \tilde{W}_0$ and its positive $\hat{\mathbf{C}}$-orbit $\hat{\omega}(y)$ can be naturally identified with the phase point and positive orbit of a standard $(\nu, k-1, r)$-flow by changing the reference coordinate system, i. e. by reducing the first component of $q_1 + q_2$ and $v_1 + v_2$ to zero. Although we will be using this equivalence, but, for the sake of brevity, we are not going to introduce an extra notation for that purpose.

The meaning of (5.4) is that there exists a small, open ball $B(z_0, \epsilon_0) \subset \mathbf{M}(\nu, k-1, r)$ of the phase space of the $(\nu, k-1, r)$-flow, such that for every $y \in \tilde{K}_0$ the positive $\hat{\mathbf{C}}$-orbit

$$\hat{\omega}(y) = \left\{\hat{y}^t = \left(\hat{q}_1^t(y), \hat{q}_2^t(y)\right) | \ t \geq 0 \right\}$$

avoids the ball $B(z_0, \epsilon_0)$.

We can now define the exponentially stable, local invariant manifolds $\gamma^s(x) \subset \tilde{W}_0$ for almost every phase point $x \in \tilde{W}_0$, similarly to (4.12), as follows:

$$\gamma^s(x) = \mathrm{CC}_x \left\{ y \in \tilde{W}_0 \big| \, d(\hat{x}^t, \hat{y}^t) \to 0 \text{ exp. fast as } t \to +\infty \right\}.$$

**Sublemma 5.6.** *The codimension-one submanifold $\tilde{K}_0$ of $\tilde{W}_0$ is transversal to the invariant manifold $\gamma^s(y) \subset \tilde{W}_0$ for every $y \in \tilde{K}_0$, whenever $\gamma^s(y)$ is a manifold containing $y$ as an interior point.*

**Proof.** This transversality immediately follows from our combinatorial assumption ($j_1 \notin Z_0(x)$ or $j_1 \in \cup_{l=1}^n Z_l(x)$ for every $x \in U_0$) by the method of the proof of Lemma 3.10. We note that this is just the point of the proof of 5.1 where we use the above mentioned combinatorial hypothesis. □

It follows now from the transversality sublemma and from the Transversal Fundamental Theorem for semi-dispersive billiards (Theorem 3.6 of [10]) that for almost every phase point $y \in \tilde{K}_0$ (with respect to any smooth measure on $y \in \tilde{K}_0$) the set $\gamma^s(y)$ is a submanifold of $\tilde{W}_0$ containing $y$ as an interior point. Therefore, the union

$$E := \bigcup_{y \in \tilde{K}_0} \gamma^s(y) \subset \tilde{W}_0$$

has positive measure in $\tilde{W}_0$. However, if the size of the open set $W_0$ is chosen small enough compared to the radius $\epsilon_0$ of the avoided ball $B(z_0, \epsilon_0)$, then the already proved ball avoiding $\hat{\omega}(y) \cap B(z_0, \epsilon_0) = \emptyset$ implies that the positive $\hat{\mathbf{C}}$-orbit $\hat{\omega}(y')$ of any $y' \in \gamma^s(y)$ avoids the shrunk ball $B(z_0, \epsilon_0/2)$, see also (4.16) and Sublemma 4.18. Thus, we obtained that in the $(\nu, k-1, r)$-dynamics the positive trajectory of every phase point $y \in E$ avoids the open ball $B(z_0, \epsilon_0/2)$, and, yet $E$ has positive measure in the phase space $\mathbf{M}(\nu, k-1, r)$. This contradicts our inductive hypothesis 4.3 postulating the ergodicity of the $(\nu, k-1, r)$-flow. Hence Lemma 5.1 follows. □

**Corollary 5.7.** *Assume the inductive Hypothesis 4.3. Then there exists a slim subset $\hat{S} \subset \mathbf{M} = \mathbf{M}(\nu, k, r)$ of the phase space of the standard $(\nu, k, r)$-flow such that every phase point $x \in \mathbf{M} \setminus \hat{S}$ either belongs to the union $\mathcal{F}$ of the codimension-one submanifolds featuring in part (2) of Lemma 4.1, or the trajectory $S^{(-\infty,\infty)}x$ contains at most one singularity and it is* **sufficient**. *Furthermore, the Chernov–Sinai Ansatz (see Condition 3.1 in [10]) holds true for the standard $(\nu, k, r)$-flow.*

**Proof.** Lemma 4.1 says that for every phase point $x \in \mathbf{M} \setminus (\hat{S}_1 \cup \mathcal{F})$ ($\hat{S}_1, \hat{S}_2, \ldots$ will always denote some slim subsets of $\mathbf{M}$) the orbit $S^{(-\infty,\infty)}x$ contains infinitely many dispersive (ball-to-ball) collisions in each time direction. Furthermore, thanks to Lemma 4.1 of [10], for every phase point $x \in \mathbf{M} \setminus (\hat{S}_2 \cup \mathcal{F})$ ($\hat{S}_2 \supset \hat{S}_1$) the trajectory $S^{(-\infty,\infty)}x$ contains at most one singularity. Then key lemmas 4.4, 3.5 and lemmas 5.1, 3.10 imply that every non-singular phase point $x \in \mathbf{M} \setminus (\hat{S}_3 \cup \mathcal{F})$ ($\hat{S}_3 \supset \hat{S}_2$) has a sufficient orbit $S^{(-\infty,\infty)}x$.

As far as the singular phase points $x \in \mathcal{SR}^+$ ($\mathcal{SR}^+$ denotes the set of all phase points $x = (Q, V^+) \in \partial \mathbf{M}$ with singular collision at time zero supplied with the outgoing velocity $V^+$) are concerned, the direct analogue of Lemma 6.1 of [16] (whose proof is quite analoguous to the proof of the original Lemma 6.1 of [16]) asserts that for a generic (with respect to the surface measure of $\mathcal{SR}^+$) singular phase point $x \in \mathcal{SR}^+$ the positive orbit $S^{(0,\infty)}x$ is non-singular and sufficient. This proves that

(i) for every phase point $x \in \mathbf{M} \setminus (\hat{S}_4 \cup \mathcal{F})$ ($\hat{S}_4 \supset \hat{S}_3$) the orbit $S^{(-\infty,\infty)}x$ contains at most one singularity, and it is sufficient;

(ii) the Chernov–Sinai Ansatz holds true.
This proves Corollary 5.7. □

**Corollary 5.8.** *Assume Hypothesis 4.3. Denote by $\mathbf{M} \setminus \mathcal{F} = \bigcup_{i=1}^{s} \Omega_i$ the decomposition of the open set $\mathbf{M} \setminus \mathcal{F}$ into its connected components, where*

$$\mathcal{F} = \left\{(q_1, q_2; v_1, v_2) \in \mathbf{M} \,\middle|\, P_{\overline{\mathcal{A}}}(v_1 - v_2) \perp z_j \text{ for some } j \leq l\right\},$$

*see part (2) of Lemma 4.1.*

*We claim that every component $\Omega_i$ ($i = 1, \ldots, s$) belongs to one ergodic component of the standard $(\nu, k, r)$-flow.*

**Remark.** We note that in the case $k = \nu$ the set $\mathcal{F}$ is empty ($l = 0$) and, therefore, $s = 1$.

**Proof.** Since the complement set $\Omega_i \setminus \hat{S}_4$ is connected and it has full measure in $\Omega_i$, the statement immediately follows from the previous corollary and the Chernov–Sinai Theorem on Local Ergodicity, i. e. Theorem 5 of [22]. □

The last outstanding task in the inductive proof of our Theorem is to connect the open components $\Omega_1, \ldots, \Omega_r$ by (positive beams of) trajectories. This will be done by the method from the closing part of the proof of Lemma 4 of [9]. Namely, it is enough to construct for every vector $z_j \in \mathbb{R}^{\nu-k}$ a piece of trajectory connecting a phase point $x^{(1)} = \left(q_1^{(1)}, q_2^{(1)}; v_1^{(1)}, v_2^{(1)}\right)$ with the property $\langle v_1^{(1)} - v_2^{(1)}; z_j \rangle > 0$ by another phase point $x^{(2)} = \left(q_1^{(2)}, q_2^{(2)}; v_1^{(2)}, v_2^{(2)}\right)$ for which $\langle v_1^{(2)} - v_2^{(2)}; z_j \rangle < 0$

This can be done, however, very easily by slightly perturbing a tangential $(1,2)$-collision for which the normal vector of impact is parallel with the vector $z_j$ and, therefore, the relative velocity of the tangentially "colliding" balls is perpendicular to $z_j$. Plainly, we can perturb this tangential collision in such a way that the perturbed collision is no longer singular, and the normal vector of impact does not change. In this way we obtain the required phase points $x^{(1)}$ and $x^{(2)}$ together with a piece of trajectory connecting them.

This completes the inductive proof of our theorem. □

## APPENDIX

### A Special Orthogonal Cylindric Billiard

In the proof of Lemma 4.21 we unfolded the non-singular orbit of (4.22) (with $|Z_i|(\nu - |Z_i|) = 0$ for $i = 0, \pm 1, \pm 2, \dots$) and thus obtained the other trajectory

$$\hat{\omega}(x_0) = \left\{ \left( \hat{q}_1^t(x_0), \hat{q}_2^t(x_0) \right) \mid t \in \mathbb{R} \right\},$$

with $\hat{q}_i^t(x_0) \in \mathbb{R}^\nu / 2 \cdot \mathbb{Z}^\nu = \mathbb{T}_2^\nu$, $i = 1, 2$. The so constructed orbit $\hat{\omega}(x_0)$ has collisions only at the cylinder

$$C = \left\{ (\hat{q}_1, \hat{q}_2) \in \mathbb{T}_2^\nu \times \mathbb{T}_2^\nu \mid d(\hat{q}_1, \hat{q}_2) = 2r \right\}$$

and at the "antipodal" cylinder

$$\overline{C} = \left\{ (\hat{q}_1, \hat{q}_2) \in \mathbb{T}_2^\nu \times \mathbb{T}_2^\nu \mid d(\hat{q}_1, -\hat{q}_2) = 2r \right\}.$$

We need to understand the ergodic properties of the arising billiard flow in the configuration space $\mathbf{Q} = \mathbb{T}_2^\nu \times \mathbb{T}_2^\nu \setminus (C \cup \overline{C})$. For the sake of simplifying the notations, we will work with the unit torus $\mathbb{T}^\nu = \mathbb{R}^\nu / \mathbb{Z}^\nu$, instead of $\mathbb{T}_2^\nu$, and we will supress ˆ from over $q$ and $v = \dot{q}$.

It turns out that the $\nu$-dimensional generator (constituent) spaces

$$\mathbb{T}_+ = \{(x, x) \mid x \in \mathbb{T}^\nu\}, \quad \mathbb{T}_- = \{(y, -y) \mid y \in \mathbb{T}^\nu\}$$

of the cylinders $C$ and $\overline{C}$ are obviously orthogonal. Thus, our system is an *orthogonal cylindric billiard* in the sense of [23]. Let us write the configuration point

$$(q_1, q_2) \in \mathbb{T}^\nu \times \mathbb{T}^\nu \setminus \left( C \cup \overline{C} \right)$$

in the form $(q_1, q_2) = (x, x) + (y, -y)$. Since the group homomorphism $\Psi : \mathbb{T}_+ \times \mathbb{T}_- \to \mathbb{T}^\nu \times \mathbb{T}^\nu$, $\Psi((x, x), (y, -y)) = (x + y, x - y)$ is surjective with the kernel $\mathrm{Ker}\Psi = \{((x, x), (x, x)) \mid 2x = 0\}$, we obtain a $2^\nu$-to-1 covering of the original flow if we switch from $(q_1, q_2) \in \mathbb{T}^\nu \times \mathbb{T}^\nu$ to $((x, x), (y, -y))$ with $q_1 = x + y$, $q_2 = x - y$. In terms of $x$ and $y$ the conditions $d(q_1, q_2) \geq 2r$ and $d(q_1, -q_2) \geq 2r$ precisely mean that $d(y, G) \geq r$ and $d(x, G) \geq r$, where $G := \{g \in \mathbb{T}^\nu \mid 2g = 0\}$. There are now two first integrals: $E_1 = \frac{1}{2}\|\dot{x}\|^2 = \frac{1}{2}\|\dot{q}_1 + \dot{q}_2\|^2$ and $E_2 = \frac{1}{2}\|\dot{y}\|^2 = \frac{1}{2}\|\dot{q}_1 - \dot{q}_2\|^2$

Thanks to the orthogonality, the $(x, \dot{x})$ and $(y, \dot{y})$ parts of the covering dynamics evolve independently of each other. This means that — after fixing the values of $\|\dot{q}_1 + \dot{q}_2\|^2$ and $\|\dot{q}_1 - \dot{q}_2\|^2$ — the $2^\nu$-to-1 covering flow is the product of two, $\nu$-dimensional, dispersive Sinai billiards. Therefore, the original orthogonal billiard flow is mixing and all but two of its Lyapunov exponents are nonzero. Besides that, for almost every phase point $z = (q_1, q_2; \dot{q}_1, \dot{q}_2)$ of this orthogonal cylindric billiard (with fixed values of $\|\dot{q}_1 + \dot{q}_2\|^2$ and $\|\dot{q}_1 - \dot{q}_2\|^2$, of course) the $(2\nu - 2)$-dimensional, exponentially contracting and expanding manifolds $\gamma^s(z)$ and $\gamma^u(z)$ exist, and they contain $z$ as an interior point.

**Acknowledgement.** A significant part of this work was done while the author enjoyed the warm hospitality and inspiring research atmosphere at the Department of Mathematics of The Pennsylvania State University, University Park Campus.